\documentclass{ceb-l}
\usepackage{amsmath, amssymb, amsthm, amsfonts, marginnote, stmaryrd,mathabx}

\usepackage{comment}

\usepackage{wrapfig}
\usepackage{graphicx}
\usepackage{xcolor}

\input xy
\xyoption{all}

\usepackage{hyperref}

\swapnumbers
\numberwithin{equation}{section}

\theoremstyle{plain}

\theoremstyle{definition}

\newtheorem{remark}[subsection]{Remark}

\def\module{\operatorname{-mod}}

\def\module{\mathrm{mod}}

\def\AA{\mathbb{A}}

\def\CC{\mathbb{C}}

\def\FF{\mathbb{F}}

\def\PP{\mathbb{P}}

\def\QQ{\mathbb{Q}}

\def\RR{\mathbb{R}}

\def\ZZ{\mathbb{Z}}

\newcommand\cA{\mathcal{A}}
\newcommand\cB{\mathcal{B}}

\newcommand\cD{\mathcal{D}}

\newcommand\cF{\mathcal{F}}

\newcommand\cH{\mathcal{H}}

\newcommand\cL{\mathcal{L}}

\newcommand\cN{\mathcal{N}}
\newcommand\cO{\mathcal{O}}
\newcommand\cP{\mathcal{P}}

\newcommand\fg{\mathfrak{g}}

\renewcommand\a\alpha
\renewcommand\b\beta
\newcommand\g\gamma
\newcommand\G\Gamma
\renewcommand\d\delta
\newcommand\D\Delta

\newcommand{\wh}[1]{\widehat{#1}}
\newcommand\quash[1]{}

\newcommand\bs{\backslash}

\newcommand{\beq}{\begin{equation}}
\newcommand{\eeq}{\end{equation}}

\newcommand{\actson}{\circlearrowright}

\usepackage[mathscr]{euscript}

\newcommand{\Gv}{{\check{G}}}
\newcommand{\Tv}{{\check{T}}}
\newcommand{\Mv}{{\check{M}}}
\newcommand{\Cx}{\CC^{\times}}

\title[The Geometric Langlands Correspondence]{What Is The Geometric Langlands Correspondence About?}

\author{David Ben-Zvi}
\address{University of Texas at Austin}
\email{benzvi@utexas.edu}
\date{\today}
\subjclass[]{}
	
\date{\today}

\dedicatory{Illustrations by Elliot Kienzle\\  https://chessapig.github.io/}

\begin{document}


\begin{abstract}
The recent proof of the unramified Geometric Langlands Conjecture has attracted a lot of publicity, so this
seems like a good time to address the title question. In one line, the Geometric Langlands correspondence is an algebraic spectral theorem for a certain class of differential equations called {\em automorphic sheaves}: it asserts they can be decomposed into monochromatic objects, which diagonalize the action of natural symmetries ({\em Hecke operators}), and it describes the corresponding colors or frequencies ({\em Langlands parameters}).

The statement is very technical and esoteric sounding, the proof takes thousands of pages, and
there are relatively few easily stated immediate consequences. So what’s the deal? In this brief survey I
will present the subject informally as a blueprint for a master plan for the study
of nonabelian symmetry, touching on some of the main motivations, connections and structures
that have emerged.

A video companion to this article is available at https://www.youtube.com/watch?v=yxauvlBZCAU
\end{abstract}

\maketitle

\tableofcontents



The spectral theorem and the Fourier transform are ubiquitous in mathematics and physics. Fundamentally they have to do with the exploitation of abelian symmetries, such as spatial translations, to decompose functions into simple, monochromatic constituents. One can imagine that a counterpart for {\em nonabelian} symmetry would be even more powerful and useful. However, our current understanding of nonabelian duality is reminiscent of the parable of the blind men and the elephant.
Number theorists have discovered a grand vision for nonabelian duality, the Langlands program, with echoes throughout representation theory and numerous applications (the most famous is the proof of Fermat's Last Theorem). Independently physicists have produced another grand vision, a nonabelian form of the electric-magnetic duality of Maxwell theory known as Montonen-Olive S-duality, with its own echoes throughout mathematical physics. Thanks to research into the Geometric Langlands Correspondence --- {\bf GLC, from now on} --- we now understand that these two communities are in fact exploring the same giant pachyderm, and
they are now able to share observations and coordinate their explorations. In fact we may now have identified the animal's beating heart -- the hidden commutativity provided by the notion of factorization algebra.
 
 This exploration has reached a major benchmark with 
the recent proof of the unramified geometric Langlands conjecture by Gaitsgory and Raskin with Arinkin, Beraldo, Campbell, Chen, Faergeman, Lin and Rozenblyum. This proof is laid out in what amounts to three papers ~\cite{GLCI,GLCIV,GLCV} and a two-volume monograph~\cite{GLCII,GLCIII}, all building on a vast body of prior technical work, including Lurie's foundations of higher algebra~\cite{HTT,HA} and much of Gaitsgory's work over the past twenty years.  To understand the statement --- let alone the proof --- of GLC requires getting through many layers of leathery hide and navigating many moving parts. However I believe the insights and perspectives gained 
amply reward the effort, and suggest a far greater whole. 

I will take this occasion to survey what kind of statement the GLC is, why we might care and why it is true. This is by no means a detailed account of the proof. My emphasis is on the parts that make the subject most appealing to me -- its multitude of connections and perspectives, and the deep underlying structures it has helped uncover. I set up the story with a discussion of the spectral theorem and Fourier transform in \S~\ref{what kind} followed by the specific context of the GLC in \S~\ref{what does it say} (and finally the statement in \S~\ref{finally GLC}). 
In \S~\ref{why true} and \S~\ref{Hitchin} I discuss two of the main structures -- factorization and the Hitchin system -- that suggest the GLC might be true, followed by an outline of the proof in \S~\ref{proof}. I conclude in \S~\ref{why} with some of the applications --- a reader impatient to see what we're doing this all for might want to glance here first --- and what might be next. 

{\bf Caveat lector:} I consistently err on the side of intuition over precision, and in references to survey articles over original sources. I have not come close to being comprehensive, or even adequate, in referencing relevant papers and important developments. I also do not discuss Phase I of the GLC Universe, from the origin of the problem in~\cite{DrinfeldGLC,LaumonGLC} to the construction of Hecke eigensheaves for $GL_n$ in~\cite{FGV,Gaitsgoryvanishing} -- see~\cite{FrenkelCEB,GaitsgoryBourbaki} for earlier surveys of the area -- but focus on the greatest hits of Phase II, the categorical GLC, from~\cite{BD} to its recent proof. I enthusiastically recommend~\cite{scholzeGLC} for an insightful in-depth account of the current state of the geometric Langlands program.

{\bf Acknowledgments:} I would like to thank the Current Events Bulletin organizers Daniel Erman and Bianca Viray. I would also like to thank Alexandre Afgoustidis, Ron Donagi, Nigel Higson, Masoud Kamgarpour, JiWoong Park and Sam Raskin for very valuable feedback on earlier drafts. I am grateful to the support of the National Science Foundation through grant DMS-2302356.



\section{What kind of statement is geometric Langlands?}\label{what kind}

I will describe the geometric Langlands correspondence as a nonabelian counterpart to the Spectral Theorem, which we now step back to revisit.

\subsection{The Spectral Theorem}
The spectral theorem provides a theory of {\em spectral decomposition} for self-adjoint operators $H$ on Hilbert spaces $\cH$. This provides a general context for the familiar fact that symmetric real matrices, or more generally hermitian complex matrices, can be diagonalized. That fact is interpreted as saying that the underlying Hilbert space $\cH=\CC^n$ decomposes as a direct sum of eigenspaces $$\cH\simeq \bigoplus_{\lambda\in \RR_t} \cH_\lambda.$$ We can think of each $\cH_\lambda$ as living over the corresponding point $\lambda\in\RR_t$, so that we can define subspaces corresponding to any subset of $\RR_t$ as the sum of the corresponding eigenspaces. 
(Here we'll denote by $\RR_t$ a copy of $\RR$ with coordinate $t$, our notation for the domain of eigenvalues and spectral measures, i.e., momentum space.) We can use this description to let all functions on $\RR_t$ act on $\cH$ by commuting operators through evaluation on the spectrum ($f$ acts by multiplication by $f(\lambda)$ on $\cH_\lambda$), so that the action of $H$ is given by multiplication by the coordinate function $t$ -- the functional calculus. 

For a general (possibly unbounded) self-adjoint operator, while $H$ might not have any actual eigenvectors in $\cH$, the spectral theorem applies measure theory to provide a sense in which we can write all vectors in $\cH$ as continuous linear combinations -- integrals over $\RR_t$ -- of eigenvectors with eigenvalue $\lambda\in \RR_t$, in effect spreading $\cH$ out as a family of vector spaces over $\RR_t$.  This can be formulated in the language of {\em projection-valued measures} on $\cH$, or more abstractly using
{\em measurable fields} of Hilbert spaces over $\RR_t$. This amounts to defining a measure $d\cH(t)$ on $\RR_t$ valued in Hilbert spaces, so that 
we recover the Hilbert space $\cH$ with the operator $H$ as the direct integral (or ``global sections'')

\begin{equation}\label{direct integral}
\cH=\int^{\oplus}_{\RR_t} d\cH(t) 
\end{equation}
while the integral on an interval $I\subset \RR_t$ produces the subspace $\cH_I \subset \cH$ where the spectrum of $H$ lies in $I$.
In this description, the operator $H$ is given by the action of the coordinate function $t$.
This action extends to the Borel {\em functional calculus}, whereby all Borel-measurable functions act on $\RR_t$ by commuting (possibly unbounded) operators -- for example the indicator function of an interval $I$ projects onto $\cH_I$. 
Variants include the continuous functional calculus, making $\cH$ a module for the commutative $C^*$-algebra $C_c(\RR_t)$, or the bounded functional calculus, a module for the commutative von Neumann algebra $L^\infty(\RR_t)$. 

An alternative ``exponentiated'' point of view on the pair $H\actson \cH$ (courtesy of the functional calculus) is as defining a unitary representation $x\mapsto e^{ixH}$ on $\cH$ of the group $\RR_x$ (a dual copy of $\RR$ with coordinate $x$ -- the position space). From this perspective, the spectral theorem provides a complete description of the unitary representations of $\RR_x$, breaking them up into atomic pieces -- the irreducible unitary representations, which are 1-dimensional eigenspaces $\CC_\lambda$. In this way the spectral theorem provides the blueprint for describing the representation theory of other groups. 

However Hilbert spaces equipped with self-adjoint operators, or with unitary representations of $\RR_x$, are living, breathing objects, and don't just form a set -- they have internal symmetry, and they can talk to each other through bounded linear operators commuting with the action (for example through unitary operators). In other words, they form a {\em category}, in fact (via the functional calculus) a category of modules over an operator algebra. While not usually phrased this way, the spectral theorem upgrades to define an {\em equivalence of categories} 
\begin{center}
\fbox{
\{Hilbert spaces with a self-adjoint operator\}$\longleftrightarrow$ \{measurable fields of Hilbert spaces on $\RR_t$\} 
}
\end{center}
-- i.e., we can spectrally decompose not just the objects but also the morphisms as linear operators between the multiplicity spaces. Again this is a model for how we might seek to describe other representation theories. 



\subsection{The spectral theorem vs. the Fourier transform}\label{symbiosis}
The spectral theorem has a symbiotic relation with the Fourier transform.
On the one hand, arguably the fundamental feature of the Fourier transform $$L^2(\RR_x)\stackrel{\sim}{\longrightarrow}L^2(\RR_t),\hskip.3in  \wh{f}(t)=\int_{\RR} f(x)e^{-2\pi i xt} dx$$
is that it takes translation, differentiation
and convolution 
operators in $x$ into multiplication operators in $t$ -- i.e.,  it simultaneously diagonalizes the action of the group $\RR$ (a unitary representation), its Lie algebra (a self-adjoint operator) and its group algebra (an operator algebra).
In the language of the spectral theorem, the Fourier transform says that the spectral decomposition of the {\em regular representation} $\RR_x\actson L^2(\RR_x)$ (or equivalently of the self-adjoint operator $i \frac{d}{dx}$) is given by $L^2(\RR_t)=\int^{\oplus}_{\RR_t}  \CC_t dt$, the ``constant rank one'' family of Hilbert spaces on $\RR_t$, with multiplicity one for every frequency. 
Said another way, Fourier inversion 
\begin{equation}\label{fourier inversion}
f(x)=\int_{\RR} \wh{f}(t)e^{2\pi i xt} dt
\end{equation} lets us build functions as superpositions of exponentials, the eigenfunctions for translation, which are the pure frequency (monochromatic) waves parametrized by frequencies (colors) $t$. 


On the other hand, the spectral theorem can itself be seen as a kind of ``categorified'' Fourier transform, in which complex-valued functions or measures and their integrals~\eqref{fourier inversion} are replaced by vector-space-valued measures and their direct integrals ~\eqref{direct integral}. In this analogy the monochromatic waves $e^{i \lambda x}$, taken by Fourier transform to $\delta$-functions $\delta_\lambda$, are replaced by the irreducible representations of $\RR_x$, i.e., the abstract eigenspaces $\CC_\lambda$, which are taken by the spectral decomposition to skyscrapers (vector-space valued $\delta$-functions) $\CC_\lambda$ at $\lambda\in \RR_t$. 
The spectral theorem then lets us 
write any unitary representation as a direct integral (superposition) of the vector space valued function given by its eigenspaces. The regular representation $L^2(\RR_x)$ plays the role of white light, a superposition of all possible frequencies, and indeed its transform $L^2(\RR_t)$ is supported everywhere on the $t$-line. 

We summarize the highlights of this analogy in the table below:

\begin{center}
\begin{tabular}[c]{l | r}
   Fourier Transform & Spectral Theorem\\
        \hline
        $e^{i\lambda x}$ & $\CC_\lambda$ \\
          $f(x)$ & $ H\actson \cH$ \\
         integral of exponentials & direct integral of eigenspaces \\
         numbers & vector spaces \\
         Hilbert space $L^2(\RR)$  & category of unitary reps \\

 \end{tabular}
\end{center}

While this perspective on the spectral theorem may appear fanciful, 
it turns out to be extremely useful. Specifically, thinking of spectral theorems as variants of the Fourier transform suggests powerful ways to characterize and then prove them -- as we will see in the case of the GLC. 

\subsection{Great, but where is Geometric Langlands?}\label{mellin}
We now massage a variant of the spectral theorem into a peculiar form, where it provides our first instance of the GLC.

\subsubsection{Spectral decomposition in algebraic geometry}
Algebra has its own form of spectral decomposition -- an elementary abelian ancestor of the GLC.
This counterpart is on the one hand less refined (and much easier) than the spectral theorem, in that we replace Hilbert spaces just with plain complex vector spaces, without any topology or inner product. On the other hand, it is more general in that we can take an arbitrary linear operator (and we could equally replace $\CC$ by any other field).\footnote{For example, an $n\times n$ nilpotent Jordan block corresponds to the $\CC[t]$-module $\CC[t]/t^{n}$, which represents a family of 1-dimensional vector spaces over the $(n-1)$st infinitesimal neighborhood of $0\in \CC_t$.
} 
The ``algebraic spectral theorem'' takes the form of an equivalence  
\begin{center}
\fbox{
\{vector spaces $V$ with an endomorphism $H$\}$\longleftrightarrow$ \{algebraic families of vector spaces over $\CC_t$\} 
}
\end{center}
Here algebraic geometry provides the notion of algebraic family of vector spaces on any algebraically defined complex variety $X$, formally known as a {\em quasicoherent sheaf} on $X$. The collection of such algebraic families, which we will encounter repeatedly, forms a category we denote $QCoh(X)$. A quasicoherent sheaf attaches a vector space (the space of sections) to the open subsets available in algebraic geometry -- those defined by the nonvanishing of polynomial functions (we don't have access to indicator functions of intervals). These vector spaces are modules for the corresponding ring of regular functions, and the theory of localization of modules tells us how to recover sections on small opens from those of larger ones. Thus algebraic spectral decomposition amounts to the polynomial form of the functional calculus: a linear operator $H\actson V$ extends uniquely to an action of the polynomial ring $\CC[t]\actson V$, whence by localization to a quasicoherent sheaf on the affine line. 

\subsubsection{The Mellin transform: Betti version} 
Let's consider the spectral decomposition of vector spaces with an {\em invertible} operator -- 
in other words, representations of the group $\ZZ$, $Rep(\ZZ)$.
The algebraic spectral theorem identifies these with modules for the ring $\CC[t,t^{-1}]$ of Laurent polynomials in one variable, i.e., quasicoherent sheaves on the punctured complex line $\AA^1_t\setminus \{0\}=\CC^\times_t$ -- we have just excluded $0$ from the spectrum. (If we throw in unitarity we get instead families of Hilbert spaces over $U(1)\subset \CC^\times_t$.) 

\begin{figure}[h!] 
        \centering 
        \includegraphics[width=0.8\linewidth]{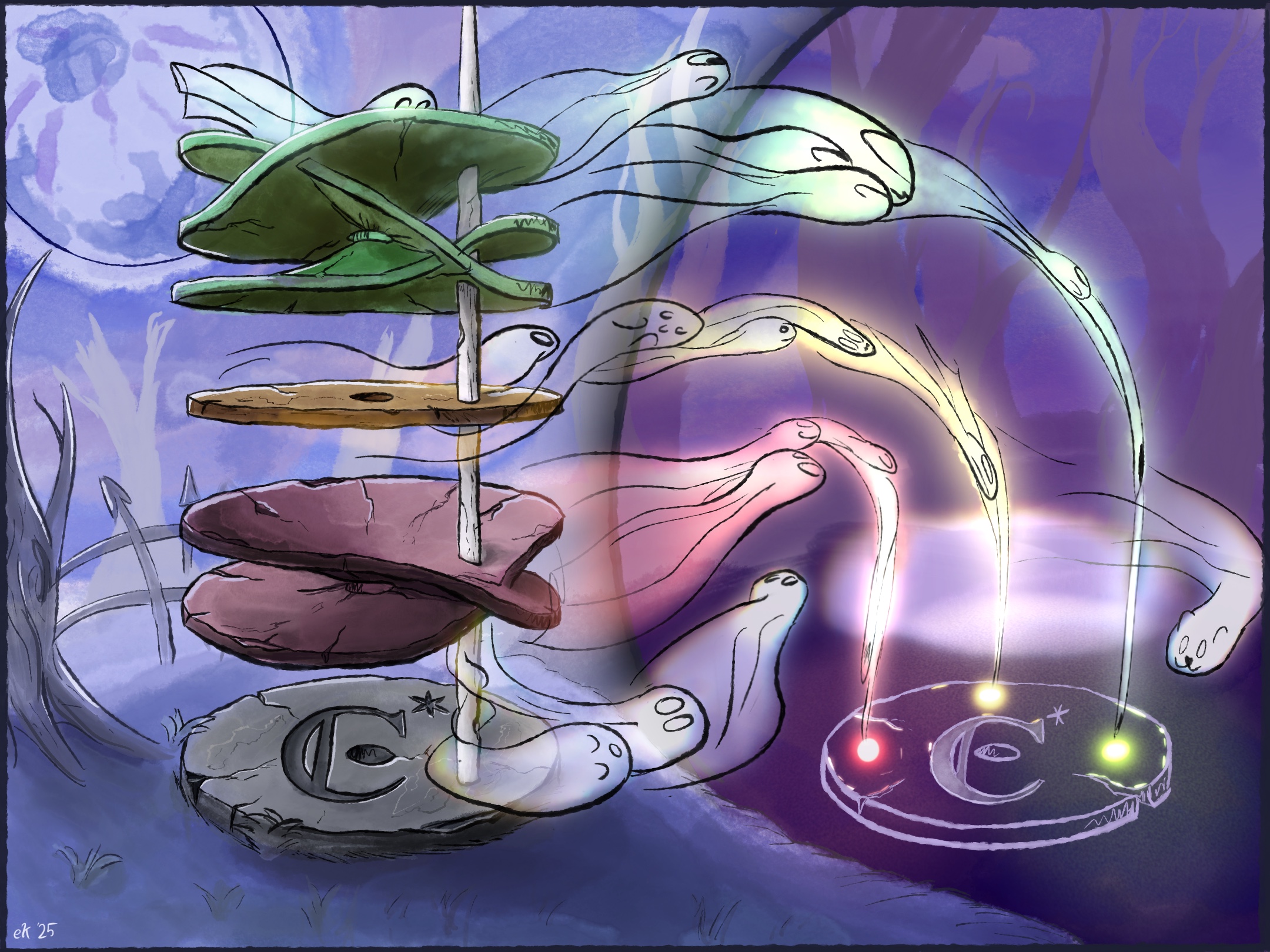} 
        \caption{The Mellin transform for sheaves: Spectres guide monochromatic local systems on $\CC^\times_z$ to their colors in $\CC^\times_t$. A skyscraper corresponds to the ``color'' white.} 
        \label{fig:Mellin} 
    \end{figure}

This result has a useful geometric interpretation, illustrated in Figure~\ref{fig:Mellin}. Let's identify $\ZZ$ with the fundamental group $\pi_1(\CC^\times_z)$ of the punctured $z$-line. Then a representation of $\ZZ$ can be described as describing the twisting or {\em monodromy} 
of a {\em local system} of vector spaces on $\CC^\times_z$. This means a locally constant family of vector spaces -- the fibers at nearby points are identified -- but following these identifications along a non-contractible loop can lead to a discrepancy, the monodromy. 
The simplest local systems on $\CC^\times_z$ are rank one families $\cL_t$ where the monodromy is multiplication 
by a nonzero scalar  $t$. These correspond to the eigenspaces (or irreducible representations) $\CC_\mu$ for $\ZZ$, parametrized by points
in the $t$-line $\CC^\times_t$. . 
The spectral theorem then gets reinterpreted as the statement that any any local system $L\in Loc(\CC^\times)$ is realized as a superposition (direct integral) of the $\cL_t$. In other words, any $L$ is encoded by the algebraic family $t\mapsto \wh{L}(t)$ of multiplicity spaces over  $\CC^\times_t$, and we 
have an equivalence, the {\em Betti Mellin transform} for local systems
\begin{center}
\fbox{
$Loc(\Cx_z)=$\{Local systems $L$ on $\CC^\times_z$\}$\longleftrightarrow QCoh(\Cx_t)=$\{quasicoherent sheaves $\wh{L}$ on $\CC^\times_t$\}  
}
\end{center}

\subsubsection{The Mellin transform: de Rham version}\label{dR Mellin}
The same spectral decomposition appears naturally in the study of differential equations. Namely, the solutions of an algebraic system of differential equations  
\begin{equation}\label{regular eqn}
z\frac{d}{dz}-A(z)
\end{equation} on $\CC^\times_z$ are typically multivalued, i.e., they undergo monodromy as we go around the circle, and so define a local system. This relation, identifying finite rank local systems with systems of differential equations with coefficients polynomial or regular functions (in our case Laurent polynomials in $z$), is an instance of the Riemann-Hilbert correspondence.
The building blocks $\cL_\mu$ correspond to the differential equations $L_s=\{z\frac{d}{dz}f=s f\}$ for a function to be homogeneous of degree $s$, solved by the multivalued function $f(z)=z^s$. (These functions are illustrated as branched covers of $\CC^\times$ in the figure.)
The multivaluedness is given by the exponential $t=\exp(s)$, which depends only on the class of $s$ in $\CC/\ZZ$. The same is true of the isomorphism class of the system $L_s$ -- conjugating by $z^n$ shifts $s$ by $n$. 

The classical Mellin transform $f(z)\mapsto \mathcal M(s)$ -- a multiplicative variant of the Fourier transform -- decomposes arbitrary functions of $z$ as superpositions of these homogeneous functions $z^s$. Likewise the building blocks $L_s$ can be used to build up arbitrary systems of differential equations with polynomial coefficients on $\CC^\times_z$. This results in a {\em de Rham Mellin transform} for differential equations:

\begin{center}
\fbox{
$\cD-mod(\Cx_z)=$\{Alg. systems of diff. eqs. on $\CC^\times_z$\}$\longleftrightarrow QCoh(\CC_s/\ZZ)=$\{quasicoherent sheaves on $\CC_s/\ZZ$\}  
}
\end{center}

While the Betti and de Rham Mellin transforms are distinct statements algebraically, they are extremely close to each other -- indeed, complex analytically the exponential map gives an identification
$$\exp:\CC_s/\ZZ \stackrel{\sim}{\longrightarrow} \CC^\times_\mu,$$
and the two versions of the Mellin transform agree on their ``common core'' of finite rank local systems and the corresponding systems \eqref{regular eqn}. The relation between the two is a variation of the relation between translation and differentiation, or representations of Lie groups and Lie algebras, mediated by the (nonalgebraic) exponential map. .

{\bf And there we have it:} the Mellin transform is itself a special case of the geometric Langlands correspondence\footnote{In the GLC, we get to fix a group --- here the group is $GL_1$ -- and a Riemann surface -- here $\Cx_z$ (or more officially, the Riemann sphere $\PP^1$ with tame ramification at $0,\infty$).}, a general spectral decomposition for certain classes of differential equations, or of their topological monodromy data, into basic monochromatic building blocks. 

What's special about these objects $\cL_t\leftrightarrow L_s$?
$\CC^\times_z$ is a group, and these objects are known as {\em character sheaves} since they behave like {\em characters} of this group -- e.g., their solutions $z^s$ are precisely characters. We can express this idea as follows: using the group structure we can translate local systems around under rotation by $w\in \CC^\times$. 
The $\cL_t=\{z^s\}$ are the {\em eigensystems}, the systems consistently isomorphic to their rotation by every $w\in \CC^\times$ thanks to $(wz)^s=w^s z^s$.
An elaboration of this statement is that the Mellin transform identifies the structure of {\em convolution} on $Loc(\CC_z^\times)$, expressing the group structure of $\CC^\times$, with the structure of {\em multiplication}, or pointwise tensor product, on $QCoh(\Cx_t)$. In this way it is both a form of the spectral theorem and a categorified Fourier transform.

This fits in a general theme of categorical Fourier transforms  -- these are equivalences of (derived) categories of sheaves of various flavors between an abelian group and a dual, exchanging convolution and tensor product operations~\cite{LaumonFourier}. Examples include the Fourier-Malgrange, Fourier-Sato and Fourier-Deligne transforms identifying $\cD$-modules, constructible sheaves and $\ell$-adic sheaves (respectively) on dual vector spaces over $\CC$, $\RR$ and $\FF_q$, and the Fourier-Mukai transform identifying coherent sheaves on dual abelian varieties~\cite{Mukai}. The challenge of geometric Langlands is to develop a nonabelian counterpart to this rich theory.

\subsection{The shape of nonabelian duality}\label{landscape}
The basic outlines of Fourier theory and the spectral theorem apply
whenever commuting operators appear. The known nonabelian counterparts are far more restrictive and delicate. 
In this section we'll sketch the roughest outline of the nonabelian dualities provided by the Langlands program
and electric-magnetic S-duality, so as to point out the privileged position enjoyed by the GLC between the two. See also \S~\ref{physics section} for a slightly broader hint as to the relation to physics. 


\subsubsection{The group.} We first fix 
 
 $\bullet$ $G$ a {\em reductive algebraic group} - for example $GL_n$, $SL_n$, $SO_{2n}$, $SO_{2n+1}$ ... $E_8$..  Over $\CC$, where we'll mostly stay, these are precisely the complexifications of compact Lie groups. (We will only consider connected groups.)

Each such $G$ comes together with a ``dual'' complex\footnote{More generally we get a group over the {\em field of definition} $k$ of the representations of $G$ we wish to study -- which we'll assume (or pretend to assume) is $\CC$.} reductive group, the {\em Langlands dual group} $\Gv$. The general idea behind $\Gv$ is that representations of $G$ of any flavor (i.e., homomorphisms {\em out} of $G$) should correspond to conjugacy classes in -- or more generally, homomorphisms {\em into} -- its dual. It can be defined concretely via the combinatorial theory of Lie groups (roots, Dynkin diagrams and so on), by keeping the Weyl group of $G$ but dualizing all of the torus and root space combinatorics. The groups $GL_n$ and $SO_{2n}$ are self-dual, simply connected groups like $SL_n$ are exchanged with adjoint groups like $PGL_n$, but the odd orthogonal groups $SO_{2n+1}$ are exchanged with the symplectic groups $Sp_{2n}$. 
As we will see, the geometric Langlands correspondence provides a fundamental description of where the dual group {\em comes from}, see \S~\ref{geometric Satake}.

\subsubsection{The setting.} Next we need to choose a setting. 

$\bullet$ Electric-magnetic duality concerns 4-dimensional quantum field theory.
For its much coarser topological form (a kind of static imprint of the theory in which we ignore all dynamics) 
we get to choose as ``spacetime'' any manifold $M$ of dimension at most 4. 

$\bullet$ For the Langlands program we get to choose one of a few flavors of field $F$: a number field like $\QQ$; a local field like $\RR$, $\CC$ or $\QQ_p$; or the field of rational functions on a (smooth projective) algebraic curve $C$. For curves $C$ over the finite field $\FF_q$, the function field $F=\FF_q(C)$ behaves very similarly to a number field. For curves $C$ over an algebraically closed field like $\overline{\FF}_q$ or $\CC$, the function field is surprisingly more akin to a local field.
 
 The GLC takes place over a smooth projective curve $C$ over $\CC$. We can equivalently think of $C$ as defining a compact Riemann surface, in particular a 2-manifold. Thus it fits neatly into both worlds. More generally there is a somewhat fanciful but extremely useful parallel between the two worlds given by the dictionary of arithmetic topology, or ``knots and primes'' analogy~\cite{MorishitaBook}, in which both number fields and curves over finite fields are thought of as (exotic relatives of) 3-manifolds, while local fields and curves over algebraically closed fields are thought of as 2-manifolds.

\subsubsection{The duality.} \label{duality intro}
To each such inputs -- a group $G$ and a field $F$ or manifold $M$ -- is associated a duality conjecture. Let us try to discern a general pattern
among all of them. On one side we have the ``question'': the automorphic or A-side, which presents a spectral decomposition problem associated to the group $G$ and the manifold $M$ or field $F$. 
These come in the form of harmonic analysis on what initially look like unrelated exotic contexts -- locally symmetric spaces $[G]_F$ associated
to number fields, moduli spaces of bundles $Bun_G(C)$ associated to algebraic curves, or moduli spaces of solutions to the Yang-Mills equations on $M$ and their variants. Remarkably these spaces are all closely related, and bridged by the GLC, as we discuss in \S\ref{BunG section}. 

On the other side we have the ``answer'': the spectral or B-side. This describes the possible colors appearing in the spectral decomposition in terms of the algebraic geometry of character varieties $Loc_\Gv$ -- spaces of representations -- of a suitable group, such as the Galois group of $F$ or the fundamental group of $M$, into the Langlands dual group $\Gv$. We discuss character varieties in the GLC setting of Riemann surfaces / algebraic curves in \S\ref{character varieties}.

Perhaps the most easily recognizable of these spectral problems is the Local Langlands Correspondence (LLC), the duality for a local field $F$ such as $\RR$ or $\QQ_p$. Here the automorphic sides consists of the collection of all {\em smooth} representations of the group $G(F)$, a real or p-adic Lie group, as opposed to the classical spectral theorem for (unitary) representations of $\RR$. The LLC seeks to describe these representations as families of vector spaces (multiplicity spaces) over the space of Langlands parameters $Loc_\Gv(F)$, playing the role of the real line $\RR_t$ in the spectral theorem. While the LLC seems immune to the much more baroque structures encountered in other settings, the field has been completely transformed in the past decade (most notably~\cite{farguesscholze}) by thinking of local fields as analogs of curves or 2-manifolds and importing ideas from the GLC and eventually from physics. See also \S~\ref{LLC}.

\subsubsection{Question and Answer?}
Labeling the two sides as question and answer conveys the psychology of the spectral theorem, and largely agrees with the flow of information in the GLC, where the geometry of character varieties of surfaces is much more elementary than the harmonic analysis of sheaves on moduli spaces. 

However in the context of number fields the typical flow of information is in the reverse direction: Galois groups and their representations are extremely mysterious and subtle objects, used to express many of the deepest questions in number theory such as solutions to Diophantine equations. Realizing them as ``colors'' of modular and automorphic forms provides a radically different description of arithmetic invariants (such as L-functions) and opens them up to an influx of powerful tools from representation theory and analysis. 

On the other hand, in physics the relation between the two sides is fundamentally different:
electric-magnetic duality in quantum field theory is a completely symmetric relation, an automorphism of a single theory much like the Fourier transform it generalizes. This automorphism acts on the various parameters of the theory such as coupling constants (exchanging weak and strong coupling), the gauge group (exchanging $G$ and $\Gv$) and the choice of topological twist -- the mechanism by which we discard dynamics and retain coarser topological information. By the time we consider the duality between topological field theories which relates directly to the GLC the A-side appears much more challenging than the B-side. However this gift from physics, the fundamental underlying balance between the two sides is still visible throughout, and percolates throughout the study of nonabelian duality.

\section{What does it say?}\label{what does it say}

We now take a closer look at the actual statement of the geometric Langlands correspondence. 

\subsection{The language: brave new harmonic and functional analysis}
We argued that the spectral theorem is concerned with a {\em category}, namely unitary representations of $\RR$ - we have not only a collection of objects, but they also talk to each other. Indeed
the basic objects of study in representation theory are categories, arising as 
representations of groups or algebras, and the goal is to describe them in a meaningful way, just as spectral theory builds vector spaces with an endomorphism out of eigenspaces. 

On the other hand, the analogy between the spectral theorem and the Fourier transform from \S~\ref{symbiosis}, in which vector spaces take the role of numbers and families of vector spaces the role of functions or measures, suggests we think of categories themselves as analogues of Hilbert spaces. Indeed objects in a category have an ``inner product'', their Hom space.
Of course it takes a lot more than a set with an inner product to make a Hilbert space --  we have addition and scalar multiplication, topology and completeness. Likewise we don't just want a category but one that has operations like direct sum and a notion of completeness. 
 In this analogy homological algebra takes the place of linear algebra over a field $k$: the well behaved class of categories in which we can perform basic operations are what we will refer to as derived categories (the technical names involve words like triangulated dg categories or $k$-linear stable $\infty$-categories). These basic examples are built by starting with abelian categories of modules for $k$-algebras and then ``completing for operations of homological algebra''. Homotopical algebra then provides counterparts to many basic constructions in functional analysis. 

Classical harmonic analysis describes the decomposition of vector spaces of functions under the action of symmetries. In recent decades, a rich categorical counterpart of harmonic analysis has emerged in which vector spaces of functions are replaced by derived categories of sheaves (i.e., different kinds of vector-space-valued functions). For example, we replace the role of functions such as $z^s$ by the differential equations they satisfy, like $z\frac{d}{dz}-s$, and transforms such as the Mellin transform by corresponding functors on categories of modules. 
The growing zoo of different sheaf theories (also known as six-functor formalisms) 
provides a counterpart to the menagerie of different function spaces. 
We summarize some of the loose analogies below:
\begin{center}
\begin{tabular}[c]{l | r}
vector spaces &  categories\\
 Hilbert spaces& derived categories\\  
operators & (exact) functors\\
 functions & sheaves\\
 Integral & pushforward / direct integral \\
 multiplication & tensor product\\
 Function spaces & derived categories of sheaves\\

 \end{tabular}
\end{center}

The proof of the GLC called for the development of counterparts for many classical chapters in functional and harmonic analysis in a brave new categorified setting. This explains some of the length and complexity of the papers in the subject, but is also one of the main dividends of the endeavor for a broader community in algebraic and analytic geometry, homotopy theory, mathematical physics and representation theory.
Due to its technical nature I will suppress as much of this language as possible in what follows.

 \subsection{The setting: moduli of bundles}\label{BunG section}.
We fix a compact Riemann surface (smooth complex algebraic curve) $C$, and a complex reductive group $G$ -- which the reader is welcomed, without much loss, to take to be the group $GL_n\CC$ of invertible matrices. 
We will focus on the {\em unramified} GLC -- meaning we don't puncture $C$ or decorate with extra data at points of $C$. This is the base case, and also the one that's been proved, though most of the general ideas have ramified versions.

The main object of study is $Bun_G(C)$, the moduli space of algebraic (or holomorphic) $G$-bundles on $C$ (denoted $Bun_G$ when $C$ is implied). This is the object of algebraic geometry which parametrizes principal $G$-bundles on $C$ (i.e., rank $n$ vector bundles for $G=GL_n$). This means it's defined by prescribing not just its points but all maps into it --- to give a map from a variety $S$ to $Bun_G(C)$ is the same as giving a $G$-bundle on $S\times C$. 
In other words, it's an algebraic geometer's version of a classifying space in topology -- specifically, it is the space of algebraic maps from $C$ to the classifying space for $G$. 

To spoil the eventual punchline, the space $Bun_G$ looks somewhat baroque from the perspective of geometry, i.e., of points or classical mechanics, but will turn out to look much nicer from the perspective of its linearization, i.e., of harmonic analysis or quantum mechanics. 

\subsubsection{The abelian case: Jacobians}
To understand $Bun_G$ we first consider the case $G=GL_1$, where its study goes back to Abel, Jacobi and Riemann. In this case a $G$-bundle is the same data as a holomorphic {\em line bundle} on $C$. Classically line bundles (or rather their holomorphic sections) are prescribed by giving a divisor on $C$, an integer linear combination of points of $C$, which determines a space of meromorphic functions with at most prescribed poles and at least prescribed zeros. Two divisors describe the same line bundle if they differ by the divisor coming from an actual rational function on $C$.

The moduli space of line bundles $Bun_{GL_1}(C)=Pic(C)$ is the Picard group of $C$. And that's already the crucial fact: the Picard group is an abelian {\em group}, under the operation of tensor product of line bundles. More concretely, the group structure is giving by addition of divisors: $Pic(C)$ is a quotient of the group of divisors, the free abelian group on points of $C$, by ``global'' relations (coming from divisors of meromorphic functions).
Moreover this group has a very simple description: it is a product $Pic(C)\simeq \ZZ\times Jac(C)$ of the integers (giving the degree of a divisor) 
by the Jacobian of $C$, which is as a topological group a torus of dimension $g$ the genus of $C$, $Jac(C)\simeq (S^1)^{2g}$.

\subsubsection{What does $Bun_G$ look like?}
\begin{figure}[h] 
        \includegraphics[width=.3\linewidth]{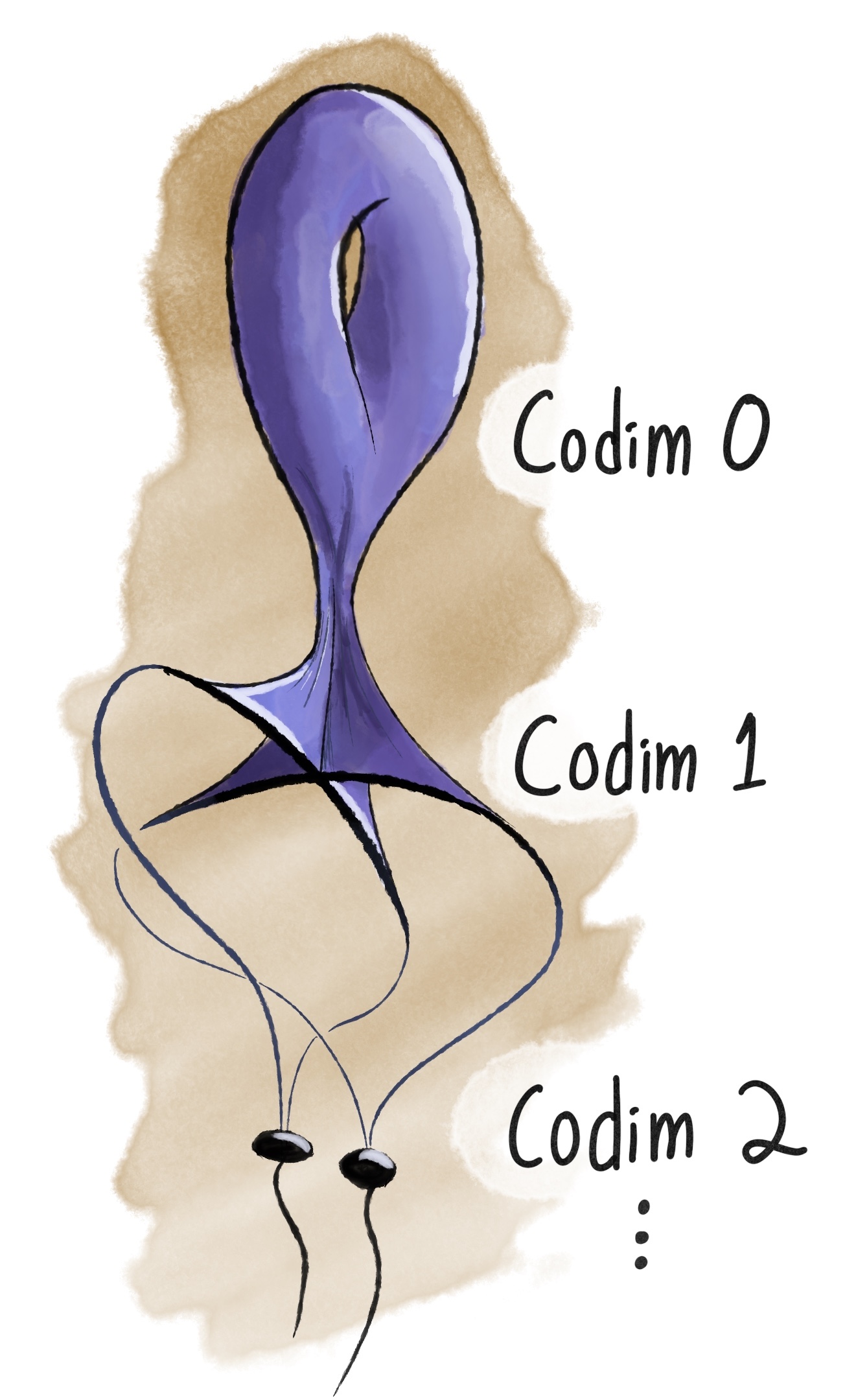} 
        \caption{The moduli of bundles and its cusps, with strata of increasing codimension.} 
        \label{fig:BunG} 
    \end{figure}

Let's now consider the geometry of the moduli spaces of $G$-bundles, which for $G=GL_n$ means the space parametrizing all holomorphic vector bundles of rank $n$ on $C$. 
A vision of the moduli space $Bun_G$ is presented in Figure~\ref{fig:BunG}. 

This space is a lot less welcoming than the Jacobian: 

$\bullet$ First, there's no natural multiplication on the collection of rank $n$ bundles (or for other $G$), so $Bun_G$ is not a group, let alone an abelian group. 

$\bullet$ Second, $Bun_G$ is not a variety but a {\em stack} -- a hybrid of classical varieties and classifying spaces, which poetically means a variety where the points carry groups of automorphisms, since bundles have potentially interesting automorphism groups. In the case of the Jacobian it is harmless to ignore these automorphisms, which are the same $\Cx$ for all line bundles (so amount just to an extra factor on a similar footing to the $\ZZ$ of connected components of $Pic(C)$).  

$\bullet$ Third, unlike the Jacobian, the connected components of $Bun_G$ are not (quasi)compact. Rather it contains an infinite ``tail'' of bundles with more and more automorphisms, built by successive extensions of line bundles over an ever-increasing range of degrees. This is closely analogous to the cusps (hence noncompactness) of modular curves such as the moduli of elliptic curves $\mathbb H/PSL_2\ZZ$.

In algebraic geometry it is common to lop off the tail of unstable bundles and then collapse the open core of semistable bundles onto a projective variety,  the moduli {\em scheme} of semistable $G$-bundles from geometric invariant theory. However for the GLC it is essential to consider the full $Bun_G$ with all its warts. First, we're doing representation theory, so discarding internal groups of symmetries is counter to the entire spirit of the subject; second, the variety of semistable bundles itself can be singular, unlike $Bun_G$, and, more crucially, has less symmetry: it is not preserved by Hecke correspondences \S~\ref{Hecke}, which are the most essential feature of the Langlands program. 
Thankfully the tails have a nice inductive structure that contrasts with the irreducible nature of the core:
the ``least compact'' part fibers over the space $Bun_T$ of bundles for the maximal torus of $G$, which is a product of Picard groups. In other words, far out on the tail bundles are given by successive extensions of line bundles. In general the tails are similarly associated with bundles for other Levi subgroups (block diagonal matrices in the case of $GL_n$ -- i.e., successive extensions of lower rank vector bundles).

\subsubsection{Divisors?}
The moduli space $Bun_G$ enjoys an important but subtle counterpart of
the concrete description of $Pic(C)$ by divisors, i.e., of sections of line bundles by functions with prescribed zeros and poles. Namely, any bundle can be trivialized away from finitely many points. Hence the bundle can be described by specifying how this trivialization fails around each of those points. This failure is captured not by an integer but by a coset for the {\em loop group} $LG=Map(D^\times, G)$ is the {\em loop group} of maps from a punctured disc into $G$ by its subgroup $LG_+=Map(D,G)$ of maps that extend to the disc. This coset space
$Gr=LG/LG_+$ is the {\em affine Grassmannian}, a counterpart of the familiar Grassmannians in the infinite-dimensional context of loop groups. Thus we can specify a bundle by specifying finitely many points on the curve labeled by points of the affine Grassmannian (the analog of a divisor), and quotienting out by the action of rational maps from $C$ to $G$ (the analog of divisors of rational functions on $C$). Thus $Bun_G$ (like the more familiar modular curves) has an ``ad\`elic'' description familiar [only] to number theorists as a double coset space $$Bun_G = G(F)\backslash G({\mathbb A}_C)/G(\cO_{{\mathbb A}_C})$$ (here $F$ is the field of rational functions on $C$). This looks unwieldy but gives more insight into the symmetries of $Bun_G$. 

For $G$ semisimple (e.g., $SL_n$ -- i.e., considering vector bundles with trivial determinant) this description has a less democratic but much more efficient variant. Fixing a point $x\in C$, $Bun_G$ can be written like a locally symmetric space in analogy with modular curves $\Gamma\backslash \mathbb{H}=
\Gamma\backslash SL_2\RR/SO(2)$, namely
$$Bun_G(C)=G(C\setminus x)\backslash LG/LG_+,$$ i.e., the quotient of a single affine Grassmannian by the group of maps into $G$ from the punctured curve. 

\subsubsection{Physics perspective}
As a brief detour, I'll mention that $Bun_G$ is closely related to the study of the Yang-Mills functional on Riemann surfaces~\cite{AtiyahBottYM}. Most famously, the moduli space of semistable $G$-bundles arises (by the Narasimhan-Seshadri theorem) as the space of minima to the Yang-Mills functional on the space of all $G$-gauge fields (connections) on $C$. However, the critical points of the Yang-Mills functional contain also unstable bundles, and the Morse theory of the functional controls the topology of the entire space $Bun_G$. Thus $Bun_G$ arises naturally in the study of quantum Yang-Mills theory on $C$. This relation is at the heart of the link between the Langlands program and gauge theory. 


\subsection{The spectral side}\label{character varieties}
On the other side of the GLC are the {\em spaces of Langlands parameters} or $\Gv$-local systems $Loc_{\Gv}(C)$ on $C$. These spaces play the role of the space of frequencies or colors $\RR_t$ in the spectral theorem, with points labeling all possible eigenvalues. We will see in \S\ref{spectral action section} how Langlands parameters appear organically as the ``answer'' from the study of harmonic analysis on $Bun_G$.

Local systems come in two main flavors, corresponding to the emphasis on translation and differentiation in Fourier theory:

$\bullet$ Betti (group) version: $Loc^B_\Gv$ is essentially the familiar {\em character variety}, parametrizing representations of the fundamental group of $C$ into $\Gv$ (i.e., of $n$-dimensional representations for $G=\Gv=GL_n$). It is a quotient (in the stacky sense) by $\Gv$ of the (mildly singular) affine variety of all homomorphisms $\pi_1(C)\to \Gv$. Concretely for $C$ of genus $g$, we can write
$$Loc_\Gv(C)=\{A_1,B_1,\dots,A_g,B_g\in \Gv\; : \; \prod_i^g [A_i,B_i]=1\}/\Gv.$$
Thus for $G=\Gv=GL_1$ this is just\footnote{Or rather, up to the stacky $\Cx$ of automorphisms, and some ``derived correction'' to account for the nontransversality of its defining equation.} $(\Cx)^{2g}$, while for genus 1 and $G=\Gv=GL_n$ this is the space of conjugacy classes of commuting matrices. The distinction with the standard notion of character variety is that, as was the case of $Bun_G$, we keep track not just of points but of automorphisms, i.e., we consider the {\em stack} of local systems rather than collapsing it to an underlying variety. 

$\bullet$ de Rham (Lie algebra) version: $Loc^{dR}_\Gv$ parametrizes $\Gv$-bundles with a flat connection, i.e., an action of the Lie algebra of vector fields, rather than the fundamental group. Thus for $G=\Gv=GL_n$ this is the moduli space of rank $n$ flat vector bundles. It forms an affine bundle (by forgetting the connection) over a quasicompact open subset of $Bun_\Gv$ (in fact it forms a twisted version of the cotangent bundle of $Bun_\Gv$). 

A flat connection defines a notion of parallel transport, its monodromy, which defines a representation of $\pi_1$. As a result, the Betti and de Rham versions of $Loc_\Gv$ are the same {\em analytically} but are distinct algebraically  -- passing to parallel transport (solving a differential equation) is a form of the exponential map. For example for $GL_1$ in genus 1 the former is $(\Cx)^2$ and the latter is a $\CC$-bundle over an elliptic curve (a non-affine surface famously considered by Serre). In other words, as far as points (or infinitesimal deformations) are concerned, the two are indistinguishable, but when considering {\em families} of local systems we have to take care.

I will default to the neutral notation $Loc_\Gv$ for both flavors of the character variety,  suppressing the distinctions whenever possible.

\subsection{Geometric Class Field Theory, or GLC for $GL_1$}
Let's (re)consider the Jacobian $Jac(C)$, i.e., the connected component of the Picard group $Pic(C)=Bun_{GL_1}(C)$. Since it is an abelian group we can perform Fourier theory on it. Specifically we can essentially repeat the Mellin transform description of \S\ref{mellin} $2g$ times. Instead of the spectral theorem for invertible operators, i.e., representations of $\ZZ$, we consider $2g$ commuting invertible operators, i.e., representations of $\ZZ^{2g}\simeq \pi_1(Jac(C))$. We can reinterpret these operators as the monodromies of a local system on the Jacobian torus $Jac(C)\sim (S^1)^{2g}$ of $C$ - the possible multivaluedness of vector-valued functions along the $2g$ directions. Spectrally, these $2g$ commuting operators as describing a module for the ring of Laurent polynomials $\CC[z_1^{\pm 1}, \dots, z_{2g}^{\pm 1}]$ in $2g$ variables, or more suggestively, an algebraic family of joint eigenspaces, a quasicoherent sheaf on $(\Cx)^{2g}$. 

On the other hand $(\Cx)^{2g}\simeq Loc_{GL_1}(C)$ is\footnote{Again, ignoring the $\Cx$ of automorphisms.} the rank one character variety of $C$ -- the data of a single scalar twisting for every loop in $C$. 
The result is an equivalence, the Betti GLC for $GL_1$, 

\begin{center}
\fbox{
$Loc(Jac(C)) = Rep(\ZZ^{2g})\stackrel{\mathbb L}{\longleftrightarrow} QCoh((\Cx)^{2g}=Loc_{GL_1}(C)).$
}
\end{center}
This equivalence has the ``Fourier'' property that translation using the group structure on the Jacobian is taken to multiplication (tensor product of modules). While our description of the equivalence uses nothing about the geometry of the curve, $\mathbb L$ satisfies a natural geometric compatibility. The inclusion $C\hookrightarrow Jac(C)$ by the Abel-Jacobi map, i.e., by considering points as divisors\footnote{This inclusion naturally lands in the degree one component $Pic^1(C)$ but we can shift it down to $Jac(C)=Pic^0(C)$ by subtracting a chosen reference point of $C$, obtaining the classical Abel-Jacobi map.}, identifies the fundamental groups (up to abelianization) $$\pi_1(C)^{ab}\simeq \pi_1(Jac(C))$$ and $\mathbb L$ sends a rank one local system $E$ on $Jac(C)$ to the skyscraper $\cO_{AJ^{-1}(E)}$ at the corresponding rank one local system on $C$. This property is the characterizing property of (unramified) {\em geometric class field theory}: $\mathbb L$ lets us identify covering spaces of $Jac(C)$ with abelian covering spaces of $C$, compatibly with the Abel-Jacobi map.

\subsubsection{$\cD$-modules and de Rham CFT}\label{dR CFT}
As with the Mellin transform there is also a de Rham version of geometric class field theory. Here instead of local systems we study {\em $\cD$-modules} on the Jacobian. Let's mention three ways of thinking of $\cD$-modules on a variety $X$, $\cD\module(X)$.:

$\bullet$ On affine varieties $X$, $\cD$-modules are simply modules for $\cD(X)$, the ring of differential operators with polynomial coefficients. When $X$ is not affine, as in the case of interest $Jac(X)$, it does not have enough polynomial functions so we use the language of sheaf theory: a $\cD$-module is a collection such modules assigned to open subsets of $X$, together with identifications on overlaps. 

$\bullet$ By thinking of presentations of $\cD$-modules by generators and relations, we can reinterpret them as {\em systems of linear PDE with polynomial coefficients}, where again we consider such systems on open patches identified on the overlaps. 

$\bullet$ A $\cD$-module carries an action of (the sheaf of) polynomial functions $\cO_X$ as well as of the Lie algebra of vector fields, interacting via the Leibniz rule. In other words, it is a quasicoherent sheaf (algebraic family of vector spaces) with a notion of differentiation, or {\em flat connection}. 

It is useful to think of $\cD$-modules as categorical avatars of generalized functions appearing as their solutions, much as we treated the equation $z\frac{d}{dz}-s$ on $\Cx$ as an avatar of its solution $z^s$. The monochromatic $\cD$-modules (character sheaves or eigensheaves) that we will be seeking are avatars of special functions and share many of their rich properties.

Among all $\cD$-modules on the Jacobian we have a special class, the line bundles with flat connection. These are in bijection with rank one local systems on $Jac$, hence (by the Betti equivalence) with rank one local systems on $C$ and finally with line bundles with flat connection on $C$. Concretely this bijection is given by simply restricting to the curve under the Abel-Jacobi map. We can interpret this bijection as attaching a $\cD$-module on $Jac$ to a skyscraper sheaf at the corresponding point on the de Rham character variety $Loc_{GL_1}^{dR}(C)$.
Laumon and Rothstein~\cite{LaumonFourier,Rothstein} proved that this attachment can be extended to an equivalence, the de Rham GLC for $GL_1$, 

\begin{center}
\fbox{
$\cD\module(Jac(C))\stackrel{\mathbb L}{\longleftrightarrow} QCoh(Loc_{GL_1}^{dR}(C))$
}
\end{center}
satisfying Fourier-like properties parallel to those of the Betti $\mathbb L$.
This equivalence is less elementary than its Betti version, in part because the situation is less affine -- neither side is presented simply by modules over a single ring. Indeed the equivalence only holds at the level of {\em derived} categories: it is an enhanced version of the quintessential derived equivalence in algebraic geometry, Mukai's famous Fourier transform for derived categories of quasicoherent sheaves on dual abelian varieties. Nonetheless it can be made very concrete and elementary for ``half'' of the parameters, which are parallel to the basic differential equations $z\frac{d}{dz}-s$ in the Mellin transform case. Namely we can write down (without any sheaf theory) explicit first order differential equations on $Jac$ by prescribing the actions of all holomorphic vector fields, $\{\xi_i f = \mu_i f\}_{i=1}^g$. These correspond to the $g$-dimensional family of all flat connections on the trivial line bundle $d-\omega$ for $\omega$ a holomorphic differential form on $C$.

\subsection{Finally, the GLC}\label{finally GLC}
We are now in a position to describe the GLC, a nonabelian generalization of the geometric class field theory from the previous section. 
The GLC studies the moduli space $Bun_G$ through linearizations, i.e., by setting up a spectral decomposition problem. The objects we will spectrally decompose -- our ``categorified function space'' -- will be particular kinds of sheaves of vector spaces on $Bun_G$. 
A sheaf on $Bun_G(C)$ can be thought of as a $C$-analog of a representation of $G$ -- indeed if we replaced $C$ by a point, the stack of $G$-bundles becomes a single point with $G$-symmetry, and sheaves of vector spaces on this are linear representations of $G$. 

So what is the kind of representation theory we will consider, the automorphic sheaves?  There are two main flavors, de Rham and Betti, parallel to the formulations of the spectral theorem in terms of self-adjoint operators (differentation) and unitary representations (translation) respectively:

$\bullet$ de Rham (Lie algebra) version: We consider $\cD$-modules on $Bun_G$, $\cD\module(Bun_G)$. As with $\cD\module(X)$ on any variety, these are defined as being compatible collections of modules for the rings $\cD$ (or systems of linear PDE) on coordinate patches, i.e., on affine varieties parametrizing families of $G$-bundles on $C$.  

$\bullet$ Betti (group) version: here we study instead {\em constructible sheaves}, families of vector spaces with the data of parallel transport along paths. The most basic examples are local systems, representations of $\pi_1$, which come from solving the parallel transport problem for flat connections (i.e. describing {\em monodromy}). But in the case of $Bun_G$ (unlike in our study of the Mellin transform or of the Jacobian) we also need to allow more singular objects where we only require parallel transport to be invertible in ``most'' directions. Such objects arise when we solve PDEs with characteristic directions that obstruct local solvability of the Cauchy problem. The precise choice of which singularities to allow (along {\em nilpotent} directions) is guided by lessons in representation theory going back to work of Harish-Chandra and Lusztig on characters of representations of Lie groups and finite groups, respectively, and will be discussed in \S~\ref{abelianization of sheaves}. The resulting category is denoted $Shv_\cN(Bun_G)$.

To avoid cumbersome notation we will refer to both categories as simply $Shv(Bun_G)$, except when we need to be more explicit about context.

Having introduced the main players, we can now state the GLC:

\medskip

{\bf Theorem.}~\cite{GLCI,GLCII,GLCIII,GLCIV,GLCV} 
There are equivalences of categories
$$Shv(Bun_G)\simeq QCoh(Loc_\Gv)$$
in both de Rham and Betti settings, as conjectured in~\cite{BD,arinkingaitsgory} and ~\cite{BettiLanglands} respectively.

\medskip

In other words, the {\em automorphic sheaves}, systems of differential equations or their topological monodromy data on $Bun_G$, have a spectral decomposition as algebraic families of vector spaces over the corresponding form of the $\Gv$-character variety of $C$. Equivalently, there is a family of special sheaves on $Bun_G$, the monochromatic objects or {\em eigensheaves}, parametrized by points of $Loc_\Gv$, and all sheaves can be written uniquely as superpositions of eigensheaves. 

Some quick notes: first, the theorem is not quite correct as stated -- we will need to modify our notion of algebraic families due to the singularities of the spaces $Loc_\Gv$, enlarging $QCoh$ to a more exotic variant called $IndCoh_\cN$, see Section~\ref{singular support}.

Second, for ``functional analytic'' purposes we have to work everywhere with {\em derived} categories -- though the equivalence is remarkably close to preserving the more familiar abelian categories of sheaves (in particular the eigensheaves themselves really are sheaves and not complexes thereof~\cite{FR}). 

Third, the theorem as stated is analogous to giving the spectral theorem as an abstract matching of two sides without characterizing it in any way (as we did with Abel-Jacobi compatibility in the abelian case) or specifying its connection to actual eigenvalues - in other words, we're still missing the main point of the whole endeavor, that this equivalence needs to {\em diagonalize something.} This is our next topic.

\section{Why is it true? The Beating Heart}\label{why true}
We have now seen the rough shape of $Bun_G$. But what is so special about it? And why should we study categorical harmonic analysis on it? Here's the key idea:

\medskip

{\bf The essence of geometric Langlands:} 
\begin{center}
\fbox{$Bun_G$ behaves surprisingly like an abelian group when linearized.}
\end{center}
\medskip

Crucially, the miraculous properties of $Bun_G$ appear not when looking at its geometry, i.e., at the level of points, but at its 
harmonic analysis or quantum mechanics, i.e., after linearization - specifically, by passing to categories of sheaves. 

This slogan has two deep manifestations that guide the subject:

\begin{enumerate} 
\item[$\bullet$]{Hecke Operators:} Sheaves on $Bun_G$ carry a large collection of commuting symmetries acting via convolution. 

\item[$\bullet$]{Hitchin fibration:} After microlocalization, i.e., passing from position space to phase space, $Bun_G$ looks generically like a Jacobian. 
\end{enumerate}

We'll start with the first of these, which explains at a fundamental level where the dual group and the duality itself come from. The second, which we pick up in  \S~\ref{Hitchin}, is much more geometric (at the cost of working only generically) and lets us see concretely what 
eigensheaves {\em look like}.


\subsection{Local symmetries: Hecke operators}\label{Hecke}
To understand the symmetries of $Bun_G$ it's useful to recall the local origin of the group structure on $Pic(C)$ as addition of divisors. Given a line bundle $\cL$ and a point $x\in C$, we can modify $\cL$ just at $x$ by considering sections that are required to vanish (or allowed to have a pole) to some order. In other words, we can add a multiple of $x$ to a divisor $\sum n_i x_i$ representing $\cL$. This gives an action of a copy of $\ZZ$ on $Pic(C)$ via $\cL\mapsto \cL(n\cdot x)$, and these actions (for varying $x$) generate the group structure.

What is the nonabelian counterpart of this $\ZZ$? Given a $G$-bundle $\cP\in Bun_G$ and a point $x\in C$ we should first parametrize all the modifications of $\cP$ at $x$, i.e., bundles that are identified with $\cP$ away from $x$. Such modifications are identified with double cosets $\underline{Gr}=LG_+\bs LG/LG_+$, i.e. by $LG_+$-orbits on the affine Grassmannian $Gr$.
A variant of Gaussian elimination shows that this set of orbits is in bijection with Weyl group orbits of characters of the dual $\Tv$ of the torus $T$ of $G$ (recovering our $\ZZ$ for $G=GL_1$). In other words, the orbits are in bijection (via the theory of highest weights) with irreducible complex representations of the Langlands dual group $\Gv$, defined combinatorially so as to have the dual torus $\Tv$ and the same Weyl group as $G$. 

But these modifications of bundles are more than just a set: we can compose modifications. Indeed while double cosets $K\backslash G/K$ in a group don't themselves carry an associative multiplication, they do so after linearization. Namely the product on the group defines an associative algebra structure on weighted combinations of group elements, i.e., measures on the group, by convolution (the group algebra). This structure survives passing to double cosets: convolution on the group defines an associative multiplication on $K$-biinvariant measures, known as a {\em Hecke algebra} for $G$ and $K$. 

We are interested in linearizing groups using sheaves rather than functions (using either $\cD$-modules or constructible sheaves). Applying this idea to the double coset space $\underline{Gr}$ gives rise to the {\em spherical Hecke category} $$Sph_G=Shv(\underline{Gr})=Shv(Gr)^{LG_+}$$
of all sheaves on double cosets for $LG_+$ on $LG$, or more concretely but less symmetrically, sheaves on the affine Grassmannian $Gr$ which are equivariant for the $LG_+$ action, or even more concretely, sheaves on $Gr$ which are constant on each $LG_+$-orbit. 
In other words, objects in the Hecke category are weighted sums of possible modifications of bundles at a point, where the weights are vector spaces. The Hecke category $Sph_G$ inherits an associative multiplication from the composition of modifications (ultimately from convolution on $LG$) -- it forms a {\em monoidal category}. Moreover, $Sph_G$ acts on $Shv(Bun_G)$ by convolution operators, the {\em Hecke functors}: we replace the value of a sheaf at $\cP\in Bun_G$ by its weighted average over bundles differing from $\cP$ by modifications at $x$. These are close relatives of the graph Laplacian, replacing the value of a function on a graph by its average over neighbors. 
 
 In the gauge theoretic description of~\cite{KapustinWitten}, taking place in the linearized setting of quantum field theory, Hecke functors are realized as {\em 't Hooft line operators} -- the modification of $\cP$ at $x$ is achieved by creating a prescribed singularity in gauge fields at $x$, which corresponds to creating a magnetic monopole passing through $x$.  
 

\subsection{Factorization: the source of commutativity}\label{factorization section}
The unmistakable beating heart of the GLC (and arguably of the entire Langlands program) is the {\bf  miraculous commutativity of Hecke modifications}. Hecke algebras in general, being averaged versions of the multiplication in nonabelian groups, have no reason to be commutative. But the physical picture of monopole operators suggests a different picture for composition, the {\em operator product expansion} (OPE): we can perform successive observations not at the same point, but at nearby points, and then take a limit as the points collide. Unlike the essential noncommutativity of quantum mechanics, there's more room to move operators around each other in higher dimensional quantum field theory. In the setting of {\em topological} quantum field theory, where measurements are constant along deformations, this makes the product of operators commutative.

The same idea is familiar in algebraic topology -- the reason the homotopy group $\pi_2$ (or $\pi_n$ for $n\geq 2$) of a topological space is abelian, unlike $\pi_1$, is that we have room to move two-fold loops in a pointed space -- maps from a disc that send the boundary to the basepoint -- around each other inside maps from a bigger disc, which makes their composition commutative up to homotopy. This observation is the origin of the theory of little disc (or $E_n$) operads (indeed of operads in general) and is the main source of commutativity in homotopy theory. 

How is this relevant? The affine Grassmannian $Gr$ is in fact homotopic to the two-fold loop space of the classifying space $BG$: it was defined as a space of $G$-bundles on a surface, and upon closer inspection describes bundles trivialized outside a small disc, hence pointed maps to $BG$. This means that $Gr$ itself carries a homotopy-commutative multiplication.

Beilinson and Drinfeld~\cite{BD} discovered a geometric unification of these mechanisms from physics and homotopy theory, known as {\em factorization} (or, synonymously, {\em fusion}), at work behind the Hecke category $Sph_G$. Namely, rather than performing successive Hecke modifications of a bundle -- acting by Hecke operators $T_1, T_2$ on sheaves on $Bun_G$ -- at the same point $x\in C$, we can instead perform the modification $T_1$ at $x$ and $T_2$ at a nearby point $y$.  The topological nature of the spherical category guarantees that the composition stays locally constant as we move $y$, or take a limit as $y$ and $x$ collide. By moving $y$ completely around $x$ we find a deformation from the composition $T_1\circ T_2$ to $T_2\circ T_1$.The result is that the convolution product on $Sph$ is actually commutative -- $Sph$ forms a {\em symmetric monoidal category}.

Factorization -- an algebraic structure parametrized purely geometrically by configurations of points on a manifold and their collisions, rather than by conventional composition maps -- was a truly revolutionary discovery, the deepest and most influential idea to come out of the study of geometric Langlands. Its origins are the efforts of Segal, Beilinson and Drinfeld to understand the geometric meaning of operator product expansion in 2d conformal field theory (as codified algebraically in the theories of vertex algebras and braided tensor categories, where the term {\em fusion} is more common), and the close parallels that emerged with the algebraic theory of iterated loop spaces (such as $Gr$) and little disc operads in algebraic topology. Once codified in~\cite{BDchiral}, it has provided a uniform perspective on phenomena in physics, topology, category theory and number theory. Turning the tables on physics, factorization is the defining feature of observables in quantum field theory in the general mathematical framework developed by Costello and Gwilliam~\cite{CostelloGwilliam}.

The commutativity of Hecke operators is arguably the most indispensable aspect of the Langlands program -- without commuting operators to diagonalize there is no spectral decomposition of automorphic forms. It is behind the construction of automorphic $L$-functions, and provides the basic mechanism linking representation theory to Galois theory. However this commutativity is proved by a trick, which fails to explain why it's there or how to turn it into a robust and malleable principle. Factorization provides precisely this principle, from which one can deduce the classical commutativity and more by decategorification. It has become a fundamental tool in arithmetic contexts thanks to the groundbreaking work of V. Lafforgue, Fargues, Scholze, Zhu and many others.

\subsection{The birth of the dual group}\label{geometric Satake}
Now that we know the spherical Hecke category is commutative, we should look for a {\em Fourier dual} description of $Sph$ in terms of tensor product rather than convolution -- i.e., a description that makes its commutativity manifest. Said more geometrically, 
we should look for its spectrum -- the space of possible (joint) eigenvalues for its actions. This is the content of the geometric Satake correspondence. 
There is already a natural guess: we noted a bijection between the double cosets 
 $\underline{Gr}=LG_+\bs LG/LG_+$ and irreducible representations of the Langlands dual reductive group $\Gv$. Thus one might hope to describe the full structure of $Sph_G$ in terms of $\Gv$-representations and their tensor product.
This guess is borne out in the following theorem:

\medskip
{\bf Theorem} [Geometric Satake Correspondence]~\cite{Lusztig, Ginzburg1995, MV}: 
There is an equivalence of symmetric monoidal abelian categories 
$Sph_G \simeq Rep(\Gv)$ between spherical category, equipped with convolution, and the representations of the Langlands dual group $\Gv$, equipped with tensor product.
\medskip

\begin{remark}
The theorem was extended in~\cite{MV} to perverse sheaves (and representations of $\Gv$) with coefficients in any field, or even over the integers. For simplicity we're implicitly passing from derived categories to their hearts, abelian $\CC$-linear categories, so that $Sph^\heartsuit$ consists equivalently of $LG_+$-equivariant $\cD$-modules or perverse sheaves on the affine Grassmanian.  There's also a vital but more subtle version of geometric Satake for the full derived category~\cite{bezfink}.
\end{remark}

The geometric Satake correspondence is perhaps better thought of in reverse: it provides the true {\em definition} of the Langlands dual group $\Gv$. Just as a commutative ring defines a geometric object, its spectrum, parametrizing all possible joint eigenvalues of the commutative ring, a commutative tensor category like $Sph$ has a spectrum given by {\em Tannakian reconstruction}. This parametrizes all possible joint eigenvalues, except these eigenvalues are now {\em vector spaces} rather than operators. 
Indeed for $G=\Gv=GL_n$ the theorem can be paraphrased as an equivalence
\begin{center}
\fbox{
\{Possible eigenvalues of $Sph_{GL_n}$\}$\longleftrightarrow$ \{$n$-dimensional vector spaces\} 
}
\end{center}
While there is only one $n$-dimensional vector space up to isomorphism, 
the theorem captures its automorphisms, which is the group $\Gv=GL_n$. In other words, the right hand side stands for the {\em stack} of $n$-dimensional vector spaces $pt/GL_n$. 
On the other hand, for $G=SO_{2n+1}$ an odd orthogonal group, with dual group $\Gv=Sp_{2n}$, the theorem says
\begin{center}
\fbox{
\{Possible eigenvalues of $Sph_{SO_{2n+1}}$\}$\longleftrightarrow$ \{$2n$-dimensional symplectic vector spaces\} 
}
\end{center}
For general $G$, the answer is given by $\Gv$-local systems on a point, which are precisely machines to turn any representation of $\Gv$ into a vector space, so a Hecke operator into an eigenvalue. 

To prove geometric Satake, we first build a group as the automorphism group of the functor of total cohomology on $Sph_G$ (which will eventually correspond to the forgetful functor $Rep(\Gv)\to Vect$). There's then a tautological tensor functor from $Sph$ to representations of this group, which we check is an equivalence.  This depends heavily on the commutativity of $Sph_G$, hence on factorization. The beautiful geometry of the Grassmannian and its Schubert stratification then forces this group to be reductive and lets us calculate its root data to be combinatorially dual to that of $G$, so we have the right to name it the Langlands dual group.

\subsection{Global symmetries: the spectral action}\label{spectral action section}
Thanks to the geometric Satake correspondence we now have a huge collection of commuting operators acting on $Shv(Bun_G)$ --- we get a copy of the commutative tensor category $Rep(\Gv)$ for every point $x\in C$, and moreover they all commute with each other. In other words, we've identified a nonabelian counterpart to the action of the group of divisors $Div(C)=\bigoplus_{x\in C} \ZZ$ on the Picard group (where the $\ZZ$ is now understood as labeling irreducible representations of $GL_1(\CC)\simeq \Cx$). That suggests we should also find a counterpart to the {\em relations} between divisors acting on the Picard, i.e., how modifications at different points interact. From the spectral point of view, geometric Satake describes the spectrum of the Hecke operators at a single point, but now we need to understand the spectrum of the ``global Hecke algebra'' generated by modifications of bundles at all points of $C$.

The key to this is again contained in the factorization structure: we have understood the composition of Hecke operators already in a way that already incorporates moving and colliding nearby points on $C$, so we only have to understand the passage from local to global. In other words, the huge collection of Hecke operators satisfy the relations that the action depends in a locally constant way on the point $x$ of modification, and performing modifications at two points that collide is the same as composing modifications at one point. This informal idea is precisely captured by the notion of {\bf factorization homology} of a factorization algebra. This idea of Beilinson-Drinfeld~\cite{BDchiral}, an abstraction of the theory of conformal blocks of vertex algebras, now provides an influential mathematical model for the way global observables in quantum field theory are assembled from local observables~\cite{CostelloGwilliam} and a powerful multi-purpose tool in topology and higher algebra~\cite{ayalafrancis,HA} as a generalization of classical homology theory.

So the challenge now, as with the geometric Satake correspondence, is to describe the globalization (factorization homology) of the spherical Hecke category in a diagonalized or Fourier-dual way, in terms of tensor product rather than convolution. In other words, we need to describe the spectrum of the global Hecke category, the collection of all its possible eigenvalues. 

Again there is a natural guess. Let's first consider the case $G=GL_n$. As we saw the possible eigenvalues for a $Sph_{GL_n}$ action at a single point $x\in C$ are given by $n$-dimensional vector spaces. The above discussion then suggests that globally we should have a vector space $\{E_x\}_{x\in C}$ which are locally constant in $x$, i.e., the vector spaces at nearby points are identified. This is precisely the data of a rank $n$ local system on $C$, i.e., a point of the character variety $Loc_{GL_n}(C)$ -- non-isomorphic local systems are distinguished by their monodromy along loops in $C$. More generally, a $\Gv$-local system $E$ on $C$ is precisely a machine that turns representations $V$ of $\Gv$ into locally constant families of vector spaces $(V)_E$ on $C$. 

This leads to the following fundamental theorem:

\medskip

{\bf Theorem} [The Spectral Action:] The Hecke action of $\bigotimes_{x\in C} Sph_x$ on $Shv(Bun_G)$ factors\footnote{Explicitly, the map $\bigotimes_{x\in C} Rep(\Gv)\longrightarrow QCoh(Loc_\Gv)$ sends a representation $V\in Rep(\Gv)$ at $x\in C$ to the vector bundle on $Loc_\Gv$ with fiber $(V)_E|_x$ at a local system $E$.} through an action of $QCoh(Loc_\Gv(C))$ in both de Rham~\cite{gaitsgoryoutline} and Betti~\cite{NadlerYunSpectral} settings.

\medskip

As with the geometric Satake correspondence, we can paraphrase this result as follows:
\begin{center}
\fbox{
\{Possible joint eigenvalues of the Hecke operators on $Shv(Bun_G)$\}$\longleftrightarrow$ \{$\Gv$-local systems on $C$\} 
}
\end{center}

The spectral action provides the {\em automorphic-to-spectral} direction of the geometric Langlands correspondence -- it tells us that the automorphic category can be spectrally decomposed over the space of Langlands parameters. Variants of the spectral action following the same general ideas have now been proved in many different settings and form some of the greatest achievements in the field, including the setting of the Langlands correspondence for function fields~\cite{LafforgueICM, AGKRRV1} and the local Langlands correspondence for nonarchimedean local fields ~\cite{farguesscholze}.  

From the point of view of gauge theory, the geometric Satake correspondence identifies the natural local observables in the theory of $G$-bundles (the Hecke = 't Hooft monopole operators) with those in a dual $\Gv$-gauge theory, built out of monodromies of $\Gv$-local systems (Wilson loop operators). The spectral action then asserts that the global observables on $C$ act on automorphic sheaves on $Bun_G$ -- in other words, it establishes a version of electric-magnetic duality on the level of observables. 

\section{How can we make it concrete?}\label{Hitchin}
In this section our narrative branches off on a (co)tangent - we briefly touch on Hitchin's integrable system, which makes parts of geometric Langlands much more geometric and explicit.

\subsection{Abelianization, or the Higgs mechanism}\label{Dolbeault GLC}

\begin{figure}[h!] 
        \centering 
        \includegraphics[width=1\linewidth]{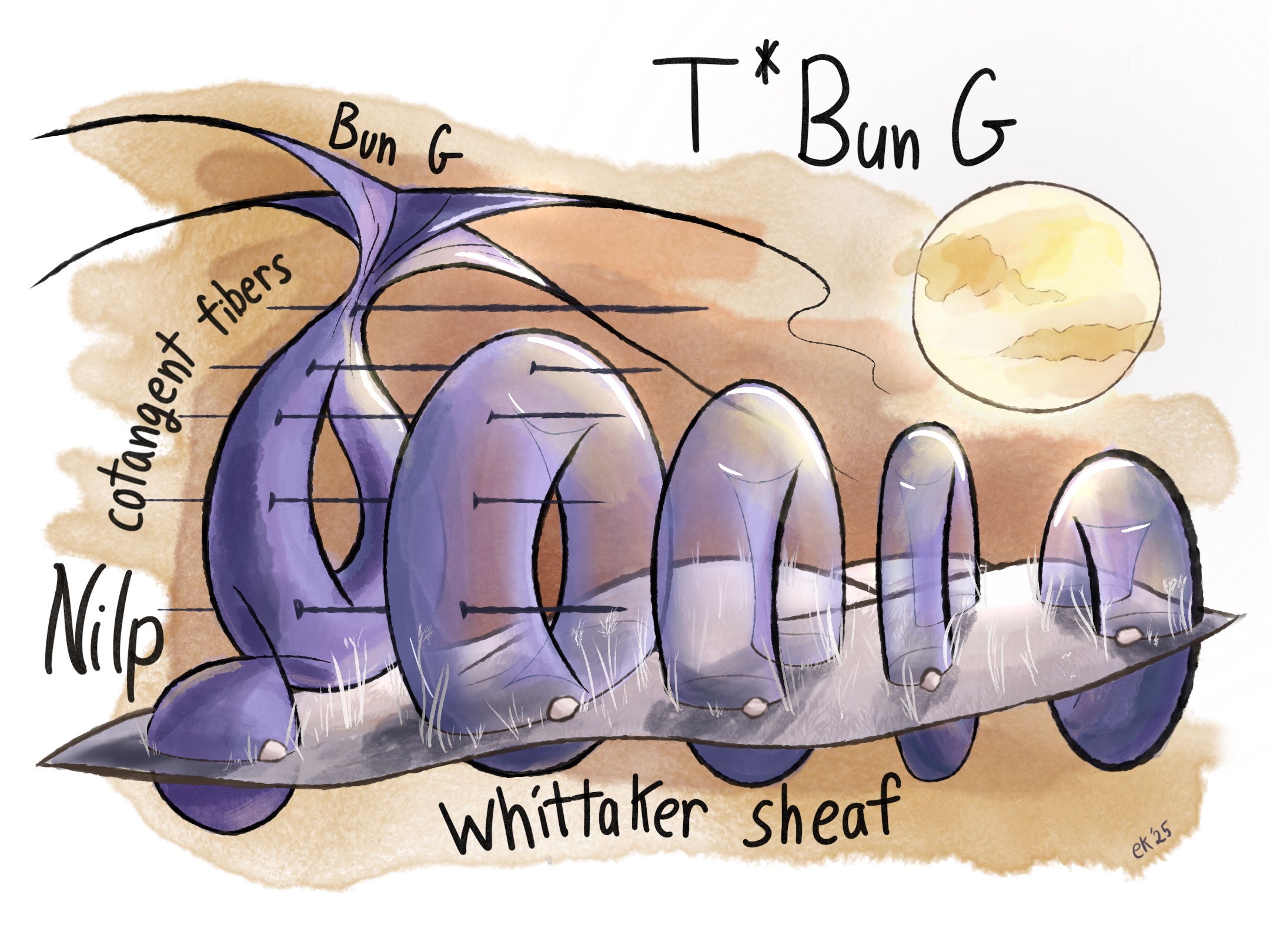} 
        \caption{The Hitchin fibration of the moduli of Higgs bundles, illustrating the zero fiber (the nilpotent cone, including $Bun_G$) and Hutchins section (microsupport of the Whittaker sheaf).} 
        \label{fig:Hitchin} 
\end{figure}

An important lesson from microlocal analysis is that the behavior of functions is closer to the geometry of {\em phase space} (where Wigner distributions, wave-front sets and characteristics of PDE live) than merely that of position space. Thus, to understand how (or even why) to perform harmonic analysis on $Bun_G$ we should really ask what's special about its {\em cotangent bundle} $T^*Bun_G$, which is illustrated in Figure~\ref{fig:Hitchin}.

The cotangent bundle of the moduli of bundles itself has a moduli space interpretation, as the moduli 
$$Higgs_G(C)=\{\cP\in Bun_G, \; \eta \in \Gamma(C, ad^*(\cP)\otimes \Omega^1_C)\}$$
 of {\em $G$-Higgs bundles}: these are bundles together with a {\em Higgs field}, an endomorphism-valued one-form (section of the coadjoint bundle twisted by 1-forms on $C$). Again we treat $Higgs_G$ as a stack, and note that it is mildly singular, like $Loc_G$. As for $Bun_G$ it also comes paired with a subquotient, the moduli space of semistable Higgs bundles, whose study was pioneered by Hitchin~\cite{HitchinSelfDuality} and is the subject of the Nonabelian Hodge Correspondence of Corlette~\cite{Corlette} and Simpson~\cite{SimpsonHiggs}. 
 
 Being a cotangent bundle, $Higgs_G$ has a natural symplectic structure (which, thanks to nonabelian Hodge theory, refines to a hyperk\"ahler structure on the corresponding semistable moduli space).  Hitchin~\cite{HitchinSystem} discovered that it also carries a gorgeous and rarefied geometric structure, that of a completely integrable Hamiltonian system, instances of which recover many of the most studied integrable systems. Let us focus for simplicity on the case $G=GL_n$, though the entire theory extends beautifully to other groups, see~\cite{donagimarkman,donagigaitsgory}. Namely, taking the characteristic polynomial of an endomorphism-valued one-form defines a map
  $$Hitch: Higgs_n\longrightarrow \bigoplus_{i=1}^{n} \Gamma(C,\Omega_C^{\otimes i})$$
from the moduli of Higgs bundles to a vector space. This map turns out to be a Lagrangian fibration, and the connected components of generic fibers $Hitch^{-1}(\chi)$ are complex tori (abelian varieties), as seen in the figure.

In fact the generic fibers are actually Picard groups (moduli of line bundles) on curves. The eigenvalues of a Higgs field $\eta$ (as a matrix-valued one-form) cut out a curve $C_\eta subset T^*C$ which is an $n$-fold cover of $C$, the {\em spectral curve}. The data of the spectral curve is equivalent to that of the characteristic polynomial $\chi(\eta)$, so we can view the map $Hitch$ as taking $\eta\to C_\eta$. The fibers then describe the possible eigen-{\em spaces} for fixed eigenvalues, which generically they a line bundle on the spectral curve -- hence the fiber is identified as the Picard group $Hitch^{-1}(\chi(\eta))=Pic(C_\eta)$. 

On the other hand the 0-fiber $\mathcal{N}ilp=Hitch^{-1}(0)$, the {\em global nilpotent cone}~\cite{LaumonNilpotent}, consists of Higgs fields where $\eta$ is a nilpotent matrix (see the figure). This is a reducible variety containing the 0-section $Bun_G\subset \mathcal{N}ilp\subset T^*Bun_G$ as the irreducible component where the Higgs field vanishes identically.

In other words, {\bf the Hitchin fibration manifests $Bun_G$ as a component of a degeneration of a Picard group of a curve} --- so $Bun_G$ is viscerally close to being an abelian variety, like it is for $GL_1$. The idea behind this abelian deformation of $Bun_G$ is a geometric implementation of a variant of the famous Higgs mechanism for mass generation and symmetry breaking in gauge theory, where the Higgs field $\eta$ (which is charged under the gauge group) acquires a vacuum expectation value $Hitch(\eta)$ and thereby breaks the gauge symmetry to the centralizer of $\eta$.

Even more strikingly, the spectral and automorphic sides of the GLC now start to look very similar: after all, connections on a vector bundle are a twisted form of endomorphism-valued 1-forms. Namely, the moduli space $Loc^{dR}_\Gv$ degenerates in a  one-parameter family to the space $Higgs_\Gv=T^*Bun_{\Gv}$ of $\Gv$-Higgs bundles. A discovery of Hausel-Thaddeus~\cite{HauselThaddeus} and Donagi-Pantev~\cite{donagipantev} is that the Hitchin systems for the Langlands dual groups generically form dual families of abelian varieties over the {\em same} vector space -- or, in the language of physics, are SYZ {\em mirrors} of each other (see the recent review~\cite{HauselICM}). Indeed ~\cite{donagipantev} proves that -- generically over the Hitchin base -- we have an equivalence, the {\em Dolbeault GLC}, 
\begin{equation}\label{DolbeaultGLC}
QCoh(Higgs_G)\longleftrightarrow QCoh(Higgs_\Gv)
\end{equation}
given fiberwise by a variant of the Fourier-Mukai transform~\cite{Mukai} identifying coherent sheaves on dual abelian varieties. (The domain of applicability of this equivalence was subsequently extended in~\cite{ArinkinFM,MRV1,MRV2}, though a full picture for the entire Hitchin system is still not available). 
This suggests that the GLC should manifest physically as a form of (homological) mirror symmetry, as was indeed confirmed in~\cite{KapustinWitten}.

\subsection{Abelianization of sheaves}\label{abelianization of sheaves}
How does the geometry of the Hitchin system relate to the theory of automorphic sheaves?

An important insight coming out of algebraic analysis is that sheaves, just like generalized functions, have a fundamentally microlocal nature, with a support in phase space. As a consequence one could approach the GLC -- as first suggested by Donagi -- by {\em abelianization}: reduce the study of sheaf theory on $Bun_G$ ``generically'' to the case of $GL_1$ using the geometry of the Hitchin fibration.

Here's the picture: in Fourier theory the monochromatic waves, the exponentials $e^{i\lambda x}$, are localized in momentum space, so look microlocally like lines $\{t=\lambda\}\subset \RR^2_{x,t}$. Likewise we might expect that the sought-after eigensheaves, the  ``monochromatic'' sheaves on $Bun_G$, should be localized microlocally along Hitchin fibers, and so accessible by abelian duality. By replacing Hitchin fibers with their limit under rescaling in the cotangent bundle, we arrive at a coarser, asymptotic version of this statement: eigensheaves should have the global nilpotent cone $\mathcal{N}ilp=Hitch^{-1}(0)$ as their characteristic variety~\cite{LaumonGLC}. (This explains to a large extent the role of sheaves with nilpotent singular support as the right generalization of local systems to consider in the formulation of the Betti GLC, \S~\ref{finally GLC}.)

In their monumental work~\cite{BD} introducing the GLC as we know it, Beilinson and Drinfeld showed that Hitchin's hamiltonians (the functions on $T^*Bun_G$ pulled back from linear functions on the Hitchin base) can be deformed to globally defined commuting differential operators\footnote{Up to a twist by spin structures.} on $Bun_G$. The eigensystems for these differential operators give a family of globally defined systems of differential equations ($\cD$-modules) on $Bun_G$ corresponding microlocally to constant sheaves on Hitchin fibers (without any winding). These are nonabelian generalizations of the differential equations $z\frac{d}{dz}f=sf$ on $\Cx$ of \S\ref{dR Mellin} and of the globally defined $\cD$-modules $\{\xi_i f = \mu_i f\}_{i=1}^g$ on the Jacobian of \S~\ref{dR CFT}.

There have since been several approaches to abelianize large parts of the GLC.

$\bullet$[Asymptotics] First, differential operators have a semiclassical limit, where they degenerate to their symbols, functions on phase space. This suggests the following idea, developed by Arinkin~\cite{Arinkinlambda}:
\begin{equation}\label{abelianization eqn}
\cD\module(Bun_G)\rightsquigarrow QCoh(Higgs_G) \longleftrightarrow QCoh(Higgs_\Gv) \leftsquigarrow QCoh(Loc_\Gv)
\end{equation}
Passing to symbols degenerates $\cD$-modules on $Bun_G$ into quasicoherent sheaves on $Higgs_G$. We then use the Dolbeault GLC ~\ref{DolbeaultGLC} to pass from $Higgs_G$ to $Higgs_\Gv$, and finally deform from $Higgs_\Gv$ to $Loc_\Gv$. This gives a ``generic, formally asmyptotic'' piece of the GLC. 

$\bullet$[Abelianization on the spectral side:] The formal WKB / semiclassical description of connections above has a striking integrated form provided by the theory of cluster varieties and from supersymmetric QFT~\cite{FockGoncharov,GMN}. This gives an abelianized description of the space of local systems $Loc_\Gv$ in terms of torus charts, coming from rank one local systems on covering space. 

$\bullet$[Reduction mod $p$] Reducing differential operators modulo a prime $p$ also has the effect of a classical limit (since $p$th derivatives vanish mod $p$, so that the ring $\cD$ becomes very close to symbols mod $p$). As a result one can relate flat connections on $C$ and $\cD$-modules on $Bun_G$ directly to their Dolbeault counterparts. This results in a proof of a generic part of the mod $p$ form of de Rham GLC by reduction to the Dolbeault case, as carried out in~\cite{BezBrav, ChenZhu}.

$\bullet$[Fukaya categories and microlocal sheaves:] Kapustin-Witten proposed that automorphic sheaves be viewed in the context of symplectic geometry and mirror symmetry as {\em A-branes}, i.e. objects of the Fukaya category, on $T^*Bun_G$. These objects are precisely represented by Lagrangians equipped with local systems. Among these are Hitchin fibers with local systems on them, providing concrete geometric realizations of Hecke eigensheaves. This had a tremendous influence on the subsequent understanding of Fukaya categories via microlocal sheaf theory~\cite{NadlerZaslow,NadlerBrane,Kontsevich,GPS}. We'll see these ideas appear in the proof of GLC, \S~\ref{Whittaker section}.

$\bullet$[Nonabelian Hodge theory:] Finally, nonabelian Hodge theory (as developed by Simpson, Mochizuki and others) leads to a profound perspective on GLC being developed by Donagi-Pantev-Simpson (see~\cite{donagipantevsurvey} for a survey as well as~\cite{DonagiPanteveigensheaves,DonagiPantevSimpson}).
The nonabelian Hodge correspondence on a variety $X$ uses harmonic map theory to identify suitable flat connections, or more general $\cD$-modules, with Higgs sheaves, or coherent sheaves on $T^*X$. Nonabelian Hodge theory on $C$ is responsible for the profound hyperk\"ahler geometry of the moduli spaces of Higgs bundles and local systems, while nonabelian Hodge theory on $Bun_G$ constructs special Hodge (or twistor) $\cD$-modules directly from line bundles on Hitchin fibers (thought of as Higgs sheaves). This provides a different way to carry out the translation of \eqref{abelianization eqn} from local systems to $\cD$-modules. See also \S~\ref{physics section}.

\section{How is it proved?}\label{proof}

In this section I'll describe some of the main ideas that go into the proof of the GLC. 

The spectral action from \S\ref{spectral action section} tells us that we can spectrally decompose the automorphic category $Shv(Bun_G)$ under the action of the Hecke operators, with the possible colors ({\em Langlands parameters}) provided by $\Gv$-local systems, i.e., by points of the character variety $Loc_\Gv$. Given this, we can formulate the ``best case scenario'':

\begin{itemize}
\item the Hecke operators have a {\em simple} joint spectrum ({\em quantum complete integrability}), and
\item every possible color appears, i.e., the spectrum is exactly the space of Langlands parameters $Loc_\Gv$. 
\end{itemize}

 In other words, we can hope that all joint eigenspaces are of multiplicity one, so that there's a {\em unique} eigensheaf $Aut_E\in Shv(Bun_G)$ (up to normalization) for {\em every} color $E\in Loc_\Gv$. The eigensheaf condition means
that every Hecke operator $\cF\in QC(Loc_\Gv)$ acts on $Aut_E$ by multiplication by the vector space $\cF_E$ (its fiber at $E$), 
$$\cF \ast Aut_E \simeq \cF|_{E} \otimes Aut_E.$$
Altogether, we can express this hope as saying that the spectral action $QCoh(Loc_\Gv)\actson Shv(Bun_G)$ is free of rank one, so that after normalization we expect to get a spectral theorem of the following form:
$$Shv(Bun_G)\simeq QCoh(Loc_\Gv).$$

Recall that a linear operator has simple spectrum precisely when it possesses a cyclic vector. 
Thus we can paraphrase this hope as saying the Hecke action has a cyclic vector (or rather cyclic object), and that this cyclic vector has full spectral support, i.e., it is a superposition of all colors. In other words, we are looking for an object $White_G\in Shv(Bun_G)$ with the properties of {\em white light}, and once we find it the GLC should be uniquely characterized by the combination
 
 $\bullet$ spectral action: the equivalence $Shv(Bun_G)\simeq QCoh(Loc_\Gv)$ respects actions of $QCoh(Loc_\Gv)$, and
  
 $\bullet$ normalization: the equivalence exchanges $White_G$ and $\cO_{Loc_\Gv}$.

\subsection{First, some functional analysis}\label{singular support}
This ``best case scenario'' turns out to need some tweaking --- i.e., we need to get the functional analysis right. 

Since the GLC is concerned with a form of harmonic analysis in which categories of sheaves replace traditional function spaces, we should encounter echos of familiar analytic issues -- and conversely, GLC provides a testing ground for this new kind of categorical analysis. As with the harmonic analysis of locally symmetric spaces, the leading issue comes from the lack of (quasi)compactness  -- i.e., what kind of growth conditions should we impose at the cusps? What classes of ``topological vector spaces'' do these choices lead to, and what kind of duality and operator theory do they carry? What do these questions even mean in the categorical setting??

 These questions were developed in the series of papers of Gaitsgory with Drinfeld and Arinkin including~\cite{DrinfeldGaitsgory, DrinfeldGaitsgorycompact, gaitsgorykernels, gaitsgorystrange, arinkingaitsgory}, based on the general categorical machinery developed by Lurie~\cite{HTT,HA}. On the one hand, they established that the basic function space in the de Rham GLC, the category $\cD\module(Bun_G)$, behaves ``miraculously'' like a categorical Hilbert space - in particular is self-dual.  
 
 On the other hand, the naive guess for the spectral side needs to be refined. In Fourier theory, lack of compact support is dual to lack of smoothness. Thus if we choose not to impose any growth conditions on our sheaves on $Bun_G$ we should expect to see singularities appearing in the sheaf theory on $Loc_\Gv$ - i.e., some kind of vector-space valued ``distributions'' rather than ``functions''. 
As it turns out this distinction on the level of sheaves is intimately related to the singularities of the underlying space $Loc_\Gv$. 

  This has led to a rich new chapter in algebraic geometry, the theory of {\em ind-coherent sheaves} $IndCoh$~\cite{GR}.
The category $IndCoh$ expands $QCoh$ to allow algebraic generalized functions. As with distributions, these new objects can be described as representing functionals on ``test'' objects, in this case, as functors on coherent sheaves. The difference between $IndCoh$ and $QCoh$ comes down to the fact that coherent sheaves on singular varieties -- for example, skyscrapers at singular points -- don't admit finite resolutions by vector bundles, or more geometrically, can't be extended to arbitrary infinitesimal deformations of the variety. 
This leads to the notion of {\em singular support} of ind-coherent sheaves, a precise counterpart of the wave-front sets of distributions, which measures the points and deformation directions where a given generalized sheaf is obstructed.

By analyzing the behavior of singular supports under natural operations on sheaves, ~\cite{arinkingaitsgory} were able to fine-tune exactly what kinds of singularities of ind-coherent sheaves match the given sheaf theory on the automorphic side: we should allow {\em nilpotent singular support}\footnote{Remarkably, this notion is not closely related to the similar-sounding nilpotent support condition imposed on Betti sheaves on the automorphic side.} -- meaning roughly that we don't allow any singularities in directions coming from tori (as in the previously understood case of $GL_1$). This finally leads to a precise statement of the unramified geometric Langlands conjecture: it's an equivalence
$$Shv(Bun_G)\simeq IndCoh_\cN(Loc_\Gv)$$ 
(in both de Rham and Betti flavors) identifying the spectral action of $QCoh(Loc_\Gv)$ on the left with its natural tensor action (multiplying distributions by functions) on the spectral side. 

\subsubsection{The \'etale setting}\label{AGKRRV section}
The landmark paper~\cite{AGKRRV1} introduced a new flavor of the GLC, the \'etale (or {\em restricted}) version. Perhaps its most important advantage is that it makes sense over arbitrary algebraically closed fields, thus formulating for the first time a categorical GLC in the setting of function fields in positive characteristic, and establishing the corresponding spectral action. Over $\CC$ the \'etale version forms an intermediate ``common core'' to the Betti and de Rham conjectures. For example on the spectral side, $Loc_\Gv^B$ and $Loc_\Gv^{dR}$ have the same points, and those points have the same deformation theory, and so the union of all the formal neighborhoods of all points is the same to both. The \'etale version of the character variety   is exactly this union for abelian groups, and a little bigger in general (essentially since exponentiating {\em unipotent} matrices is an algebraic operation). In particular coherent sheaves built in a finite process out of skyscrapers are the same on all three versions. Dually, the \'etale automorphic category is designed so as to be built {\em discretely} rather than continuously out of monochromatic objects.

The paper~\cite{AGKRRV1} proved that the Betti and de Rham versions of GLC are equivalent, by showing they both reduce to the common \'etale version. This lets us play off advantages of the two settings -- for example the Betti character stack $Loc_\Gv^{B}$ is very simple, as the quotient of an affine variety by $\Gv$, and manifestly depends only on the topology of $C$, while the de Rham version $Loc_\Gv^{dR}$ is far from affine, so that Hom spaces between coherent sheaves are finite dimensional.

\subsection{Interlude: Relative Langlands, or, Where have all the $L$-functions gone?}\label{relative GLC}
The Langlands philosophy doesn't just predict a relation between automorphic and spectral sides diagonalizing Hecke operators: it comes equipped with a huge repertoire of predictions for which measurements should match which. The study of these measurements is known as the relative Langlands program, and was systematized in~\cite{SV}. The most famous -- the matching of {\em L-functions} associated to automorphic and Galois representations -- is arguably the main feature of the arithmetic Langlands program. $L$-functions themselves are conspicuously absent from most of the literature on the GLC, but similar matching measurements play a major role in normalizing, characterizing and proving the GLC.

What do we mean by a measurement of automorphic forms? These are linear functionals on the vector space of automorphic forms --- though we can (and will) think of them dually, as represented by inner product with a particular vector, or in the geometric context, as Hom spaces from a particular sheaf. The special measurements we can make predictions for -- the {\em periods}\footnote{Confusingly, not to be confused with periods in algebraic geometry.} -- are highly structured: they have counterparts in all the different settings of the Langlands program (local, global, arithmetic and geometric) satisfying strong compatibilities. 

In physics, such structured collections of measurements or states form an essential part of any quantum field theory, where they correspond to {\em boundary theories} (or codimension one defects). More suggestively, we can study gauge theories through their interaction with charged materials. In the context of electric-magnetic S-duality their study was pioneered by Gaiotto-Witten~\cite{GaiottoWittenboundary,GaiottoWittenSduality}, and proposed as a general framework for the relative Langlands program in~\cite{BZSV} (see also~\cite{ganBZSV}).

A basic example is the {\em trivial period} -- the constant or zero-frequency wave, which corresponds to the simplest measurement: taking the average of an automorphic form or global cohomology of an automorphic sheaf. This period is already very interesting in the context of GLC. It was studied in in~\cite{gaitsgorylurie} and used to prove Weil's Tamagawa Number conjecture for function fields. 

In the next three sections, we will look at three important types of measurements, or periods, whose interplay provides the backbone of the proof of GLC (to paraphrase~\cite{GLCIII}). 
Here we list the names of the matching automorphic and spectral measurements and their primary role: 

\begin{enumerate}

\item[$\bullet$] Eisenstein $\leftrightarrow$ Eisenstein: provides an inductive structure on the GLC, reducing the conjecture to the {\em cuspidal} case.

\item[$\bullet$] Whittaker $\leftrightarrow$ global sections (Neumann): used for {\em normalization} - provides a faithful measurement of cuspidal automorphic sheaves.

\item[$\bullet$] Kac-Moody $\leftrightarrow$ opers: Used for the {\em construction} of cuspidal automorphic sheaves.
\end{enumerate}

Perhaps the most crucial feature of periods is that they satisfy a strong local-to-global compatibility -- they are given as {\em $\Theta$-series}. The $\Theta$-series (or Poincar\'e series) are generalizations of the classical Riemann $\theta$-function appearing for example in Riemann's integral formula for $\zeta(s)$.
They are modular and automorphic forms that are written as averages of local data over a global symmetry group. This in fact provides one of the most common way to write down automorphic forms. 

A characteristic feature of L-functions is their ``local-to-global'' description as Euler products. 
In the geometric setting, periods have a loosely analogous key property -- {\bf factorization, again}. Namely we can attach $\Theta$-series to local data at finitely many points of the curve $C$, and the construction is compatible with collisions of points. This factorization property is the glue that lets us assemble local inputs into a comprehensible global output. It is also technically very challenging -- a major reason for the length of the proof of GLC (and many of its predecessors) is the need to understand holistically local, global and factorizable versions of all the constructions.

\subsection{The Inductive Step: Eisenstein series} 
Eisenstein series are simultaneously the hard part and the easy part of the theory of automorphic forms.
On the one hand, the serious need for analysis comes from the continuous spectrum of
modular curves (and more general arithmetic locally symmetric spaces), an effect of their lack of compactness -- i.e., from functions not vanishing at infinity (i.e., at the cusps). The general analytic theory of Eisenstein series, one of Langlands' tour-de-force achievements, provides a construction of this continuous spectrum. The geometric Langlands counterpart of Eisenstein series has to grapple with issues parallel to the analytic difficulties in the classical theory but in the brave new world of categorical analysis. 

On the other hand, representation theoretically Eisenstein series provides the ``unsurprising'' part of the theory -- the part coming by {\em parabolic induction}. This is an elaboration of highest-weight theory which lets us build representations of big groups from those of smaller {\em Levi} subgroups (for $GL_n$ these are the subgroups of diagonal or block-diagonal matrices). The guiding light of representation theory of reductive groups, 
Harish-Chandra's philosophy of cusp forms, says that for each group we should first sift out the families (or {\em series}) of representations coming by parabolic induction from smaller groups, leaving only the rare gems -- the cuspidal representations -- which we then carry with us on our journey to ever larger groups. 

The entire Langlands program follows this inductive strategy: it dictates that the construction of Eisenstein series from Levi subgroups $M\subset G$ should match with a similar parabolic induction for the dual Levi subgroup $\Mv\subset \Gv$, which precisely accounts for the locus of {\em reducible} Langlands parameters. Thus very roughly we need to first ``subtract'' the contributions of Eisenstein series from both sides of the correspondence, leaving behind the cuspidal automorphic functions (those with vanishing Eisenstein periods), which we need to then match with {\em irreducible} Langlands parameters. 

To make a long and technical story criminally short, the theory of Eisenstein series is now very well understood in the geometric setting along these general patterns, following decades of contributions by Drinfeld, Laumon, Braverman, Gaitsgory, Finkelberg, Mirkovic, Bezrukavnikov and others. The monograph~\cite{GLCIII} completes the inductive step, showing in particular the compatibility of Eisenstein series with the Langlands functor $\mathbb L$ as above. With this in hand, the proof of GLC reduces to the study of the orthogonal complement of Eisenstein series, the cuspidal automorphic sheaves.

\subsection{White Light: the Whittaker period}\label{Whittaker section}
At the opposite extreme from monochromatic waves we find white light, a superposition of all colors. 
In the setting of the Fourier transform, white light is represented by the $\delta$-function $\delta_0$ at $x=0$, whose Fourier transform is the constant function $1$. This choice plays an essential role in normalizing eigenvectors for differentiation (or translation), since the eigenproperty characterizes them only up to scale. The normalized exponentials are chosen to have value $1$ at $x=0$, i.e., by their pairing with the delta-function $\delta_0$, so as to be constituents of white light. 

This normalization extends readily to more general abelian dualities -- we normalize eigenfunctions by evaluating at the unit. In other words we stipulate that the duality sends the delta-function $\delta_e$ at the unit (a Dirichlet boundary condition) to the constant function $1$ (a Neumann boundary condition). This is illustrated (in the case of the Mellin transform) in Figure~\ref{fig:Mellin} by the white pole corresponding to the white field. 

In nonabelian duality, experience in representation theory suggests that the role of white light is played by a less familiar object, the {\em Whittaker\footnote{Named after the mathematician Sir Edmund Whittaker. Serendipitously, the Old English etymology of the name Whittaker means White Field.} period} (known as the {\em Nahm pole boundary condition} in physics~\cite{GaiottoWittenSduality}). In other words, we seek to normalize the Langlands correspondence by stipulating that the Whittaker period is taken to the spectral measurement given by total integral or global sections (i.e., to the Neumann boundary condition). In other words, the sought-after duality should take the representing object, the {\em Whittaker sheaf}, to the structure sheaf $\cO_{Loc_\Gv}$ representing spectral global sections. 

What is this Whittaker period? It derives from Whittaker functions, solutions of the Whittaker differential equation from mathematical physics, which were used to give {\em normalized} descriptions of representations of groups over local fields. In the setting of modular forms the Whittaker period corresponds to the measurement of the first Fourier coefficient in the $q$-expansion, which is used to normalize cusp forms.
 
 In the setting of the GLC, the Whittaker period has an appealing microlocal geometric description, illustrated in Figure~\ref{fig:Hitchin} -- it is represented by the {\em Hitchin section}~\cite{HitchinSection}. This is a section of the Hitchin integrable system consisting of those Higgs fields which are in rational cyclic form as matrices for $G=GL_n$ (or more generally lie in Kostant's reductive-group generalization of rational cyclic form). These correspond to spectral curves equipped with a trivial line bundle, i.e., the identity elements of the Hitchin fibers. A key difference from the abelian case, also visible in the figure, is that the Hitchin section completely misses the zero-section $Bun_G\subset T^*Bun_G$ (since the zero matrix doesn't admit a cyclic vector). 
 
 Note that being a section this is transverse to the Hitchin fibers, which provided our microlocal picture of what eigensheaves look like -- as one might hope from a candidate for white light. 
 
 Once we've picked out the Whittaker sheaf, we can start to act on it by Hecke functors. (Microlocally, this means we move the Hitchin section around by translation along the Hitchin fibers). Thanks to the spectral action, this assignment determines a functor $$QCoh(Loc_\Gv)\to Shv(Bun_G),\hskip.3in \cF \mapsto \cF\ast Whit_G$$ sending the structure sheaf $\cO_{Loc_\Gv}$ to the Whittaker sheaf. If we really expected the automorphic side to be free of rank one over the spectral action, this {\em would} be the GLC. 
 
However we know this functor can't be quite right -- we need to fix the spectral side to the slightly bigger category $IndCoh_\cN(Loc_\Gv)$. This requires some careful analysis. The papers~\cite{FR,NadlerTaylor} used the microlocal interpretation to show some surprisingly nice properties of Whittaker coefficients (the measurement represented by the Whittaker sheaf) including {\em exactness} on a large subcategory of sheaves.
 These results are combined with further estimates in~\cite{GLCI} (concerning cohomological {\em boundedness} and categorical {\em compactness})
to show that testing automorphic sheaves against ones of the form $\cF\ast Whit_G$ ($\cF\in QC(Loc_\Gv)$) does lift uniquely to a functor
in the opposite direction $$\xymatrix{Shv(Bun_G)\ar[rr]^-{\mathbb L}\ar[rrd]&&
IndCoh_\cN(Loc_\Gv)\ar[d]\\ && QC(Loc_{\Gv}). }$$

The microlocal geometry of the Whittaker period is also used in~\cite{FR} to demonstrate a key injectivity property of the geometric Langlands correspondence which is runs somewhat contrary to experience in number theory: Whittaker coefficients don't vanish on cuspidal sheaves. Namely, we can test any cuspidal sheaf faithfully against the Whittaker sheaf and its Hecke translates, so it won't be annihilated by the functor $\mathbb L$. Combined with the inductive understanding of the Eisenstein part, what remains is to somehow establish surjectivity -- to see that every color appears. 

\subsection{Construction: Kac-Moody localization and opers}
Finally, we need to {\em find} some rare gems -- cuspidal automorphic sheaves. A lot of them, in fact, to account for all irreducible Langlands parameters. As we mentioned, the main technique to actually write down classical automorphic forms is the theory of $\Theta$-series -- i.e., to start with local data and average.
The ace up the sleeve of the de Rham geometric Langlands discovered by Beilinson and Drinfeld is a particular type of geometric $\Theta$-series, the {\em Kac-Moody period}, which has a concrete dual description in terms of {\em opers}. This construction differs from all the topics we have discussed in that it does not have any known direct counterpart in other parts of the Langlands program. 

Why is it so elusive? The answer is that while the entire de Rham GLC is based on representation theory of Lie algebras rather than Lie groups, the Kac-Moody period in particular concerns representations that {\em do not} integrate to the group, such as the action of $\fg$ on its enveloping algebra $U\fg$ itself. This kind of representation has no counterpart in the Betti or \'etale GLC or indeed in the smooth representation theory of p-adic groups, the topic of the local Langlands correspondence\footnote{There {\em are} however counterparts in the rapidly developing world of locally analytic (p-adic) representations of p-adic groups.}.

Another form of the answer comes from physics: the Kac-Moody period is {\em not topological} -- it is sensitive to the complex structure on our Riemann surface $C$. As such it is related to the much richer study of conformal or holomorphic quantum field theory rather than the relatively well understood topological field theory~\cite{FrenkelGaiotto}.

So what {\em is} this Kac-Moody period about? The global version is the sheaf of differential operators $\cD_{Bun_G}$ itself -- a very natural object but one that has no analog in the topological world of monodromies and constructible sheaves. More precisely, we need to use the correction of differential operators involving spin structures, a universal feature of quantization. The local version is the category of representations of the loop Lie algebra $\fg((z))$ -- the Lie algebra of $LG$ -- again up to a ``spin'' correction, giving what are known as {\em critical level} representations of a central extension of $\fg((z))$, the affine Kac-Moody algebra $\wh{\fg}$.
 The local-to-global construction ($\Theta$-series) is given by one of the greatest hits of geometric representation theory, Beilinson-Bernstein localization construction: a method of writing systems of differential equations prescribed by representations of a Lie algebra. Through this construction, representation theory provides a wealth of fairly concrete $\cD$-modules on $Bun_G$, which Beilinson and Drinfeld realized was the secret recipe for building a Langlands correspondence.

What does representation theory have to offer here? The starter for this recipe is a miraculous theorem of Feigin and Frenkel~\cite{FeiginFrenkel} (see also the book ~\cite{FrenkelLanglands}), using the approach to loop algebras suggested by conformal field theory, the theory of vertex algebras. This theorem identifies the intrinsic parameters for critical level representations of $\wh{\fg}$ with a special class of Langlands parameters ($\Gv$-local systems) on the punctured disc called {\em opers}. Opers (in the case of $GL_n$) are connections with a choice of cyclic vector for differentiation -- i.e., connections which are in rational cyclic form. These are exactly 
the first order systems of $n$ differential equations in one variable we get when we want to solve a single $n$-th order ODE. 
Opers form a dual, spectral counterpart of the Whittaker period -- a version of Hitchin's section for the space $Loc_\Gv$ of flat connections. 

The Feigin-Frenkel theorem provides something remarkable -- a direct bridge between representation theory associated to $G$ and the Langlands dual side, which we can use to populate the automorphic side with eigensheaves.
As we mentioned in \S~\ref{abelianization of sheaves}, Beilinson-Drinfeld write a family of explicit differential equations on $Bun_G$ given by global differential operators -- i.e., built directly out of $\cD_{Bun_G}$ -- which are nonabelian analogs of the equations $z\frac{d}{dz}f=sf$ on $\Cx$ of \S\ref{dR Mellin} and their counterparts for Jacobians \S~\ref{dR CFT}. By identifying these $\cD$-modules as the output of the local-to-global construction from those $\wh{\fg}$-representations corresponding to $\Gv$-opers on the entire curve $C$, they were able to show that they were indeed Hecke eigensheaves with precisely these opers (thought of as points of $Loc_\Gv$) as their eigenvalues. 

But this is just the beginning: opers themselves only form a half-dimensional subvariety of $Loc_\Gv$. However {\em any} irreducible local system admits the structure of an oper away from finitely many points~\cite{arinkinoper} (for $GL_n$ this is a classical result of Deligne that flat connections admit meromorphic cyclic vectors). Moreover up to homotopy there's only one such choice -- the space of choices is contractible (a theorem of~\cite{BKS} for classical groups, and a collolary of the GLC in general). So one has to ask more out of the representation theory: a description of a broader class of representations associated to $\Gv$-local systems with a rational oper structure, provided by ~\cite{FrenkelGaitsgoryspherical}. Now we're cooking. 

In~\cite{GLCII} the critical level representation theory is built into a highly structured machine incorporating factorization and local-to-global compatibility. This machine is used to populate the automorphic side with D-modules coming from Kac-Moody localization, and the spectral side with matching sheaves coming from local systems with oper structure.This gives a rich enough repertory of matching objects to see the ``surjectivity'' of the Langlands functor. 

Putting all these ingredients together we have the five paper sequence:
\begin{enumerate}
\item[\cite{GLCI}] the construction of the Langlands functor,
\item[\cite{GLCII}] the construction of cuspidal sheaves by Kac-Moody localization,
\item[\cite{GLCIII}] the inductive step provided by Eisenstein series,
\item[\cite{GLCIV}] more subtle functional analysis of the Langlands functor (it satisfies the categorical self-duality property of {\em ambidexterity}), 
\item[\cite{GLCV}] leading to the conclusion that all eigenspaces have multiplicity one and the unramified GLC holds.  
\end{enumerate}

\section{What is this good for, and what's next?}\label{why}
Thus after a herculean effort the unramified geometric Langlands correspondence is now a theorem. But as I've tried to emphasize the point of the endeavor was not solely or even primarily verifying the truth of a particular statement. Research around the GLC has led to a deeper understanding of nonabelian duality as it appears in representation theory, arithmetic and quantum field theory and to a rather surprising joint vision behind these far-flung areas. The long sequence of papers culminating in the proof doesn't try to follow a geodesic path, but rather surrounds the statement and besieges it before eventually breaking in. The proof is instead a major benchmark of our understanding of representation theory and harmonic analysis, the first time a nonabelian duality statement this exhaustive and general has been proved. It demonstrates that as a community we have understood the overarching features of nonabelian duality, including the key subtleties and roadblocks, and have an adequate conceptual framework, centered around factorization, to carry this understanding to new contexts. Indeed the main ingredients in the proof - with the critical exception of the magic of opers and Kac-Moody representations - are essentially portable.

\subsubsection{Ramified GLC}
Before moving further afield, I want to briefly discuss the meaning of the modifier ``unramified''. The colors appearing in our discussion thus far are given by local systems (de Rham or Betti) on a smooth projective curve / compact Riemann surface $C$. It is very natural to consider more generally ramified local systems, meaning bundles with {\em meromorphic} flat connections, or representations of the fundamental group of  punctured surfaces $C\setminus \{x_1, \dots, x_n\}$. These arise as eigenvalues for Hecke actions defined only away from the marked points, such as are found on moduli spaces of bundles on $C$ equipped with extra structure at the markings. Bezrukavnikov~\cite{romaicm,romahecke} proved a landmark theorem extending the geometric Satake correspondence to the setting of {\em tame} ramification (parabolic structures on bundles / regular singularities on connections). This theorem provides in particular a precise formulation of a tamely ramified version of the GLC, which currently remains open in general.

To formulate (let alone prove) a full ramified extension of the GLC requires a greater understanding of the {\em local} GLC, which provides a conjectural spectral decomposition for the categorical representation theory of loop groups in terms of Langlands parameters given by local systems on the punctured disc. See~\cite{GaitsgoryICM} for a discussion of the rapidly changing landscape around this conjecture. 

\subsubsection{Representation Theory of Reductive Groups}\label{rep theory motivation}
The GLC has a central role in representation theory, where it serves as a prototype for a spectral theory for representations associated to a reductive group $G$. What makes this field remarkable is that the same $G$ has many distinct avatars, which nonetheless share the same structural features. At the core of the subject since the 1990s is a fundamental triangle of representation theories, linked by a web of conjectures of Lusztig~\cite{CarterLusztig}: 
\begin{enumerate}
\item quantum groups $U_q \fg$,
\item loop groups $LG$ and affine Kac-Moody algebras $\widehat{\fg}$, and 
\item reductive groups and Lie algebras in positive characteristic $G, \fg/\FF_q$
\end{enumerate}
Ideas from the GLC --- specifically the geometric Satake correspondence~\cite{MV} and its tamely ramified refinement (Bezrukavnikov's theorem~\cite{romaicm,romahecke}) --- have completely transformed this area, for a very partial sample see~\cite{BMR,BezMirk,RicheGeordielinkage,AcharRiche,richesurvey,Geordiesurvey,FrenkelGaitsgorylocalization, RaskinYang} (see also \S~\ref{FLE}).

\subsubsection{GLC in low genus} When we restrict to curves of genus zero and one, both sides of the correspondence become very explicit and classical objects, for which the Langlands dual description provided by geometric Langlands provides deep insights. For example, in genus one the Langlands parameters are commuting pairs of matrices, whose algebraic geometry has long been elusive~\cite{LiNadlerYun}. On marked or punctured curves of genus zero Langlands parameters become meromorphic ODE or, in the \'etale setting, Galois representations of the field of rational functions in one variable. In this (ramified) setting, the GLC has been particularly important in the study of the special class of {\em rigid local systems}~\cite{zhiweirigidity}, including the arithmetic of Kloosterman's exponential sums and Bessel equations~\cite{FrenkelGross,HeinlochNgoYun,XinwenBessel} and combinatorics of Schubert calculus~\cite{LamTemplier}. A striking application is Yun's solution of the inverse Galois problem for some finite reductive groups~\cite{zhiweiGalois}. In the other direction, the automorphic spectral problems studied by the GLC are identified in low genus with the solution of famous quantum integrable systems, such as the Gaudin spin chain and the Calogero-Moser particle systems. The GLC then provides a powerful paradigm for solving quantum integrable systems in which the joint spectrum of commuting Hamiltonians is given by Langlands parameters~\cite{FrenkelGaudin}.

\subsubsection{Other GLCs}\label{FLE}
The GLC has both a {\em degeneration}, the Dolbeault GLC, and a {\em deformation}, the quantum GLC, which remain open even in the unramified setting, though as discussed in \S~\ref{Dolbeault GLC} the Dolbeault conjecture was proven {\em generically} in~\cite{donagipantev}.
Both have the beautiful feature that they are completely symmetric in $G$ and $\Gv$, and both appear naturally from the same physics origin (S-duality for 4d $\cN=4$ super-Yang-Mills) as the GLC.  See~\cite{PadurariuToda} for progress on a new form of the Dolbeault conjecture over the entire Hitchin base inspired by Donaldson-Thomas and geometric invariant theory.
The quantum GLC provides a global geometric context for the (noncritical) representation theory of affine Kac-Moody algebras (de Rham) and quantum groups (Betti). Its symmetric treatment of $G$ and $\Gv$ has influenced the development of the ``usual'' GLC in several ways, going back to the proof of the Feigin-Frenkel theorem from which the quantum GLC evolved. For example, the Kac-Moody and Whittaker periods, used for the ``surjectivity'' and  ``injectivity'' in the proof of GLC, have quantum variants which are on completely equal footing and in fact are exchanged  by the duality. The essential local ingredient in the quantum GLC replacing the geometric Satake correspondence is the so-called Fundamental Local Equivalence~\cite{gaitsgoryWhittaker, chenfu, FLEtrio,gaitsgorylysenkoFLE} relating representations of affine Kac-Moody algebras and quantum groups, the most refined in a sequence of developments going back to Drinfeld and Kohno's construction of quantum groups from monodromy of the Knizhnik-Zamolodchikov equations of conformal field theory and the Kazhdan-Lusztig equivalence~\cite{KazhdanLusztig}.

\subsection{Arithmetic Applications}
The geometric Langlands program was originally conceived by Drinfeld as a counterpart to the classical Langlands correspondence where one could benefit from powerful tools of algebraic geometry and sheaf theory with the hope of gaining insight into the original problems in arithmetic. 
In the past decade this hope has been realized perhaps more than could reasonably have been expected. Below I focus on two of the most successful and active arenas for this interaction, the settings of local fields and of function fields in positive characteristic. These are by no means exhaustive. I refer the reader to the excellent overview~\cite{zhuICM} for some highlights of the steady stream of recent applications of GLC to a wide range of problems in both the local and global Langlands program and arithmetic geometry more broadly, such as the Tate and Beilinson-Bloch-Kato conjectures. I also refer to~\cite{ganBZSV} for an overview of relative Langlands duality (touched on in \S~\ref{relative GLC}), where ideas from geometric Langlands and its physics interpretation have been applied to the theory of L-functions and their representations as integrals (periods) of automorphic forms. 

\subsubsection{Local Langlands}\label{LLC}
The local Langlands correspondence (previously encountered in \S~\ref{duality intro}) provides a conjectural spectral theory for representations of reductive groups over local fields $F$, like $\QQ_p$, $\RR$ and $\FF_q((t))$.
In other words it proposes an identification of representations with families of vector spaces (sheaves of multiplicity spaces) over a ``spectrum''. Moreover it identifies the parameter space or spectrum as a space of Langlands parameters -- collections of conjugacy classes in $\Gv(\CC)$, organized as representations of a form of the Galois group of $F$ into $\Gv(\CC)$. 
The LLC is a crucial partner to the global Langlands program, much as the spectral theorem is to the Fourier transform: global Langlands concerns the spectral decomposition of the particular representation $L^2([G]_F)$ of the $G(F_v)$ (for $v$ a place of a global field $F$) associated to $G$ and $F$ (for instance, $L^2(SL_2(\ZZ)\bs SL_2(\RR))$). 

Both the LLC and GLC take the form of spectral theorems describing categories (of representations or sheaves) in terms of geometry of Langlands parameters. This aligns well with the parallel provided by arithmetic topology between local fields and 2-manifolds. Indeed the study of the LLC has recently been transformed under the influence of the GLC into a categorical statement, describing {\em categories} of representations of groups over local fields as categories of quasicoherent sheaves over Langlands parameters --- this is the topic of Volume II of the recent IHES proceedings on the Langlands program~\cite{IHESLanglands}. 

This parallel has become much more direct thanks to the seminal works of Fargues-Scholze~\cite{farguesscholze} in the p-adic case and Scholze~\cite{scholzetwistor} in the real case. In these works the local Langlands correspondence is literally realized as a case of the (unramified!) geometric Langlands correspondence, but in a much more exotic setting than curves over a field  -- the Fargues-Fontaine curve and twistor $\PP^1$, respectively. Many of the ideas developed in the GLC -- most notably factorization and the spectral action -- have already played a crucial role here, notably in establishing the automorphic-to-Galois direction of local Langlands~\cite{farguesscholze} and more recently the independence of $\ell$ in the LLC~\cite{scholzemotivic}. This provides a direct path for the ideas and techniques of the GLC (and indirectly, structures from physics) to impact number theory. This is an extremely active area, with new frontiers including an emerging geometric formulation of the Langlands correspondence for p-adic analytic representations of p-adic groups. 

\subsubsection{Langlands for Function Fields}\label{function fields intro}
A very different (and much older) relation between geometric and classical Langlands comes from the close analogy (Weil's Rosetta Stone) between number fields and function fields of curves over a finite field $\FF_q$, and between Riemann surfaces and curves over its algebraic closure\footnote{From the arithmetic topology perspective, this corresponds to thinking of a 2-manifold as the fiber of a fibered 3-manifold.} $\overline{\FF}_q$. 
In this analogy vector spaces of automorphic forms arise as decategorifications of categories of automorphic sheaves, via the mechanism of {\em traces of Frobenius} (Grothendieck's function-sheaf correspondence). This mechanism is a direct extension of the link provided by the Weil conjectures between point counts for varieties over finite fields with cohomology over their algebraic closure. Under this mechanism, the spectral decomposition of automorphic forms (the classical Langlands correspondence) is a shadow of the categorical decomposition of automorphic sheaves (geometric Langlands). 

The original dream behind developing the GLC is that one might translate one's way back across Weil's Rosetta Stone, carrying along insights as we pass from Riemann surfaces to curves over finite fields and maybe eventually to number fields. Indeed the development of the GLC has largely taken place in parallel over $\CC$ and over $\overline{\FF}_q$.
In this way geometric Langlands ideas were crucial to advances such as Ng\^o's proof of the Fundamental Lemma~\cite{nadlerCEB,HalesBourbaki}, one of the major breakthroughs in the arithmetic Langlands program;  V. Lafforgue's work establishing one direction of the Langlands correspondence for function fields~\cite{CaraianiCEB,StrohBourbaki}; the proof by Gaitsgory-Lurie of Weil's Tamagawa number conjecture for function fields~\cite{gaitsgorylurie}, and the proof by Gaitsgory of de Jong's conjecture~\cite{deJongConj} on monodromy over finite fields. 

The \'etale version of the GLC of~\cite{AGKRRV1} discussed in \S~\ref{AGKRRV section} was proven in~\cite{AGKRRV3} to imply, through a sophisticated form of taking traces of Frobenius, a significantly strengthened form of the unramified Langlands conjecture for function fields. Now that unramified GLC in characteristic zero is understood, the pendulum has swung back and many of the big structural questions for function fields feel within sight. Very recently~\cite{gaitsgoryraskin2025} a large part of the unramified positive characteristic conjecture has been proved. The ICM addresses~\cite{GaitsgoryICM} and~\cite{RaskinICM} outline visions for understanding the fully ramified Langlands correspondence over function fields and for proving some of the most important open problems in the theory of automorphic forms in positive characteristic, the Ramanujan and Arthur conjectures.

\subsection{Relations to physics}\label{physics section}
As I have emphasized, a lot of what gives the GLC its appeal and power is its proximity not only to arithmetic but to physics, which I will quickly recap here. 

Kapustin and Witten~\cite{KapustinWitten} (see also~\cite{Frenkelgauge,BettiLanglands}) explained how the GLC naturally arises as an aspect of Montonen-Olive S-duality -- a conjectural nonabelian extension of the electric-magnetic symmetry of Maxwell theory.
S-duality is a non-perturbative equivalence between two four-dimensional quantum field theories, the maximally supersymmetric ($\cN=4$) extensions of Yang-Mills theory for two Langlands dual Lie groups. To get to geometric Langlands, we throw away more and more information from this incredibly rich statement. First we pass to a {\em twist} of the theory, meaning that we throw away all of the metric dependence of the theory (e.g., dynamics) by passing to cohomology of a suitable differential, provided by the structure of supersymmetry. (To get a sense for the severity of this process, its analog for Maxwell theory (the case $G=GL_1$) amounts to forgetting everything about electromagnetism except for Gauss' law!) There is a family of differentials depending on a parameter $\Psi\in \CC \mathbb P^1$, and S-duality gives rise to an equivalence  $\cA_G\simeq \cB_\Gv$ between two 4d TQFTs, an ``automorphic'' theory (the 4d A-model, $\Psi=0$)) associated to $G$ and a ``spectral'' theory (the 4d B-model, $\Psi=\infty$) associated to $\Gv$. 

We then simplify further by only considering the theory on 4-manifolds of the form $C\times \Sigma$ for a fixed Riemann\footnote{Different aspects of the construction depend either on a complex structure on $C$ or just on the underlying 2-manifold.} surface $C$. In other words, we extract from the original 4d gauge theory an assignment from Riemann surfaces to 2d topological quantum field theories. Finally, we follow the example of homological mirror symmetry and extended TQFT and focus, for each such 2d TQFT, only on its category of boundary conditions. Thus we end up with an assignment from Riemann surfaces to categories (a higher analog of the extremely rich and well-studied notion of 2d conformal field theory).

Kapustin and Witten discovered that this subsequent simplification applied to the physical S-duality conjecture produces exactly the kinds of categories appearing in the geometric Langlands correspondence, while varying $\Psi$ leads to the quantum GLC (the Betti form of GLC arose as an attempt to make this relation to TQFT explicit). Their work, and two decades of subsequent work, shows that as we pay attention to more and more of the forgotten structure from physics we can see all of the key features of the geometric Langlands program emerge as well as remarkable new features. 

First, the physics derivation from 4d gauge theory (and even further, from 6d QFT) provides a unified framework with which to understand numerous connections (many of which predate~\cite{KapustinWitten}) with 2d conformal field theory, classical and quantum integrable systems, mirror symmetry, gauge-theoretic invariants in topology and the Seiberg-Witten geometry of supersymmetric gauge theory.

Second, we find that geometric Langlands is just (part of) the 2d part of a statement about 4d theories, so extends (conjecturally) to statements involving manifolds of every dimension up to four (including with boundaries, corners and other defects) which are very tightly controlled by the algebraic structure of TQFT. The B-model TQFT $\cB_{\Gv}$ is fairly well understood, while the A-side is quite mysterious. However, thanks to the proof of GLC we now know remarkable TQFT properties of the categories of (Betti) automorphic sheaves $\cA_G(C)=Shv_\cN(Bun_G(C))$: for example, they are acted on by the mapping class group of $C$ and satisfy simple cut-and-paste rules. 
Most strikingly, it follows from the proof of GLC that there is a well-defined vector space $\cA_G(M)$ of ``automorphic forms on a 3-manifold $M$'', which can be calculated in terms of categories of automorphic sheaves on surfaces using techniques such as Heegaard splittings. It also satisfies a Langlands duality identifying it with (a stacky, derived enhancement of) the space of functions on the character variety $Loc_\Gv(M)$ (space of representations of $\pi_1(M)$ into $\Gv$). We know very little about these spaces directly, but after deforming to the {\em quantum} GLC these spaces become the much-studied skein modules of 3-manifolds associated to quantum groups, see the review~\cite{jordansurvey}, resulting in Langlands duality predictions for skein modules. 

Another place where the categorical complexity of geometric Langlands gives way to concrete isomorphisms of vector spaces is the analytic Langlands correspondence~\cite{EFK}, a remarkable theory of spectral decomposition for quantum integrable systems which arises physically (following~\cite{GaiottoWittenAnalytic}) as an instance of the equivalence $\cA_G\simeq \cB_\Gv$ on a special 3d configuration.

Third, physics demonstrates the full symmetry between the two sides of the correspondence -- after all, we start with {\em the same} 4d QFT applied to two Langlands dual groups, a symmetry that's broken by choices made in the process but can be restored by varying those choices, but can then be use to inform our understanding of both sides. This symmetry is not at all obvious in number theory, where modular forms and Galois representations seem to come from different universes, but suggests new perspectives (for example in the setting of the relative Langlands program \S~\ref{relative GLC}). 

Finally, the physics of supersymmetric QFT is much richer than that of topological field theory. In particular 
it lets us distinguish a class of objects, the {\em BPS branes}, which preserve more supersymmetry than just the differential used to define the topological theory. This richness is reflected mathematically in
the rigid structures of Hodge theory. On a Riemann surface $C$, this results in the enhancement of the moduli of Higgs bundles and local systems from holomorphic symplectic to hyperk\"ahler geometry. On $Bun_G$, it results in the enhancement of the sheaf theory from $\cD$-modules to Hodge or twistor modules. (These structures play a key role in the work of Donagi-Pantev-Simpson encountered in \S~\ref{abelianization of sheaves}.) This suggests a link to the theory of motives. Motivic sheaves, like these supersymmetric branes, have many different realizations (de Rham, Betti, Dolbeault, \'etale). They have recently come to prominence in the GLC and LLC~\cite{scholzemotivic}. More broadly, understanding motives themselves (as enhancements of Galois representations) is one of the main motivations for the entire Langlands program. It's tempting to dream that the physics of supersymmetry and S-duality might provide some insight about these profound questions.

\bibliographystyle{alphaurl}
\bibliography{CEB}

@article {deJongConj,
    AUTHOR = {Gaitsgory, D.},
     TITLE = {On de {J}ong's conjecture},
   JOURNAL = {Israel J. Math.},
  FJOURNAL = {Israel Journal of Mathematics},
    VOLUME = {157},
      YEAR = {2007},
     PAGES = {155--191},
}

@book {MorishitaBook,
    AUTHOR = {Morishita, Masanori},
     TITLE = {Knots and primes},
    SERIES = {Universitext},
      NOTE = {An introduction to arithmetic topology},
 PUBLISHER = {Springer, London},
      YEAR = {2012},
     PAGES = {xii+191}
     }

@article {Lusztig, 
    AUTHOR = {Lusztig, George},
     TITLE = {Singularities, character formulas, and a {$q$}-analog of
              weight multiplicities},
 BOOKTITLE = {Analysis and topology on singular spaces, {II}, {III}
              ({L}uminy, 1981)},
    SERIES = {Ast\'{e}risque},
    VOLUME = {101},
     PAGES = {208--229},
 PUBLISHER = {Soc. Math. France, Paris},
      YEAR = {1983},
   MRCLASS = {17B10 (05A30 20G05 22E47)},
  MRNUMBER = {737932},
MRREVIEWER = {James E. Humphreys},
}

@inproceedings {LafforgueICM,
     AUTHOR = {Lafforgue, Vincent},
     TITLE = {Shtukas for reductive groups and {L}anglands correspondence
              for function fields},
 BOOKTITLE = {Proceedings of the {I}nternational {C}ongress of
              {M}athematicians---{R}io de {J}aneiro 2018. {V}ol. {I}.
              {P}lenary lectures},
     PAGES = {635--668},
 PUBLISHER = {World Sci. Publ., Hackensack, NJ},
      YEAR = {2018},
   MRCLASS = {11R39 (11F70 14G35 14H60)},
  MRNUMBER = {3966741},
MRREVIEWER = {Shrenik Shah},
}

@article {Hecke,
    AUTHOR = {Hecke, Erich},
     TITLE = {\"{U}ber {M}odulfunktionen und die {D}irichletschen {R}eihen mit
              {E}ulerscher {P}roduktentwicklung. {I}},
   JOURNAL = {Math. Ann.},
  FJOURNAL = {Mathematische Annalen},
    VOLUME = {114},
      YEAR = {1937},
    NUMBER = {1},
     PAGES = {1--28},
      ISSN = {0025-5831},
   MRCLASS = {DML},
  MRNUMBER = {1513122},
       DOI = {10.1007/BF01594160},
       URL = {https://doi.org/10.1007/BF01594160},
}

@article {SV,
AUTHOR= {Sakellaridis, Yiannis AND Venkatesh, Akshay},
TITLE= {Periods and harmonic analysis on spherical varieties},
   JOURNAL = {Ast\'erisque},
  FJOURNAL = {Ast\'erisque},
    NUMBER = {396},
      YEAR = {2017},
     PAGES = {360},
      ISSN = {0303-1179},
      ISBN = {978-2-85629-871-8},
}

@article {NadlerYunSpectral,
    AUTHOR = {Nadler, David and Yun, Zhiwei},
     TITLE = {Spectral action in {B}etti geometric {L}anglands},
   JOURNAL = {Israel J. Math.},
  FJOURNAL = {Israel Journal of Mathematics},
    VOLUME = {232},
      YEAR = {2019},
    NUMBER = {1},
     PAGES = {299--349},
      ISSN = {0021-2172},
   MRCLASS = {14D24 (22E57)},
  MRNUMBER = {3990944},
MRREVIEWER = {Alessandro Ruzzi},
       DOI = {10.1007/s11856-019-1871-9},
}

@article {gaitsgoryoutline,
    AUTHOR = {Gaitsgory, Dennis},
     TITLE = {Outline of the proof of the geometric {L}anglands conjecture
              for {$GL_2$}},
   JOURNAL = {Ast\'{e}risque},
  FJOURNAL = {Ast\'{e}risque},
    NUMBER = {370},
      YEAR = {2015},
     PAGES = {1--112}
}

@article {gaitsgorykernels,
    AUTHOR = {Gaitsgory, Dennis},
     TITLE = {Functors given by kernels, adjunctions and duality},
   JOURNAL = {J. Algebraic Geom.},
  FJOURNAL = {Journal of Algebraic Geometry},
    VOLUME = {25},
      YEAR = {2016},
    NUMBER = {3},
     PAGES = {461--548},
       DOI = {10.1090/jag/654},
}

@article {gaitsgorystrange,
    AUTHOR = {Gaitsgory, Dennis},
     TITLE = {A ``strange'' functional equation for {E}isenstein series and
              miraculous duality on the moduli stack of bundles},
   JOURNAL = {Ann. Sci. \'{E}c. Norm. Sup\'{e}r. (4)},
  FJOURNAL = {Annales Scientifiques de l'\'{E}cole Normale Sup\'{e}rieure. Quatri\`eme
              S\'{e}rie},
    VOLUME = {50},
      YEAR = {2017},
    NUMBER = {5},
     PAGES = {1123--1162},
       DOI = {10.24033/asens.2341},
}

@book{gaitsgorylurie,
 ISBN = {9780691182148},
 author = {Dennis Gaitsgory and Jacob Lurie},
 publisher = {Princeton University Press},
 title = {Weil's Conjecture for Function Fields: Volume I (AMS-199)},
 year = {2019}
}

@article {HitchinSystem,
    AUTHOR = {Hitchin, Nigel},
     TITLE = {Stable bundles and integrable systems},
   JOURNAL = {Duke Math. J.},
  FJOURNAL = {Duke Mathematical Journal},
    VOLUME = {54},
      YEAR = {1987},
    NUMBER = {1},
     PAGES = {91--114}
}

@article {HitchinSelfDuality,
    AUTHOR = {Hitchin, N. J.},
     TITLE = {The self-duality equations on a {R}iemann surface},
   JOURNAL = {Proc. London Math. Soc. (3)},
  FJOURNAL = {Proceedings of the London Mathematical Society. Third Series},
    VOLUME = {55},
      YEAR = {1987},
    NUMBER = {1},
     PAGES = {59--126}
}

@article {SimpsonHiggs,
    AUTHOR = {Simpson, Carlos T.},
     TITLE = {Higgs bundles and local systems},
   JOURNAL = {Inst. Hautes \'Etudes Sci. Publ. Math.},
  FJOURNAL = {Institut des Hautes \'Etudes Scientifiques. Publications
              Math\'ematiques},
    NUMBER = {75},
      YEAR = {1992},
     PAGES = {5--95},
      ISSN = {0073-8301,1618-1913},
   MRCLASS = {32G13 (14D07 53C07 58D27 58E15)},
  MRNUMBER = {1179076},
MRREVIEWER = {William\ Goldman},
       URL = {http://www.numdam.org/item?id=PMIHES_1992__75__5_0},
}

@article {Corlette,
    AUTHOR = {Corlette, Kevin},
     TITLE = {Flat {$G$}-bundles with canonical metrics},
   JOURNAL = {J. Differential Geom.},
  FJOURNAL = {Journal of Differential Geometry},
    VOLUME = {28},
      YEAR = {1988},
    NUMBER = {3},
     PAGES = {361--382},
      ISSN = {0022-040X,1945-743X},
   MRCLASS = {58E20 (32L99 53C10)},
  MRNUMBER = {965220},
MRREVIEWER = {John\ C.\ Wood},
       URL = {http://projecteuclid.org.ezproxy.lib.utexas.edu/euclid.jdg/1214442469},
}

@article{BD,
AUTHOR={Alexander Beilinson and Vladimir Drinfeld},
title = {Quantization
of {H}itchin {H}amiltonians and {H}ecke Eigensheaves},
URL={https://math.uchicago.edu/~drinfeld/langlands/QuantizationHitchin.pdf}
}

@book {BDchiral,
    AUTHOR = {Beilinson, Alexander and Drinfeld, Vladimir},
     TITLE = {Chiral algebras},
    SERIES = {American Mathematical Society Colloquium Publications},
    VOLUME = {51},
 PUBLISHER = {American Mathematical Society, Providence, RI},
      YEAR = {2004},
     PAGES = {vi+375},
       DOI = {10.1090/coll/051},
}

@article {DrinfeldGaitsgorycompact,
    AUTHOR = {Drinfeld, Vladimir and Gaitsgory, Dennis},
     TITLE = {Compact generation of the category of {D}-modules on the stack
              of {$G$}-bundles on a curve},
   JOURNAL = {Camb. J. Math.},
  FJOURNAL = {Cambridge Journal of Mathematics},
    VOLUME = {3},
      YEAR = {2015},
    NUMBER = {1-2},
     PAGES = {19--125},
       DOI = {10.4310/CJM.2015.v3.n1.a2},
}

@article{DrinfeldGaitsgory,
 author = {Vladimir Drinfeld and Dennis Gaitsgory},
 title = {On some finiteness questions for algebraic stacks}, 
  journal ={Geom. Funct. Anal},
  year={2013}, 
  issue={23(1)}, 
  pages={149--294},
       DOI = {10.1007/s00039-012-0204-5},
}

@article{scholzeGLC,
      title={Geometric Langlands (after Gaitsgory, Raskin, ...)}, 
      author={Peter Scholze},
      year={2026},
      \JOURNAL = {Bourbaki Report},
              Note = {\url{https://people.mpim-bonn.mpg.de/scholze/Exp1252_Scholze.pdf}},

}

@article{farguesscholze,
      title={Geometrization of the local {L}anglands correspondence}, 
      author={Laurent Fargues and Peter Scholze},
      year={2021},
      eprint={2102.13459},
      archivePrefix={arXiv}
}

@incollection {romaICM,
    AUTHOR = {Bezrukavnikov, Roman},
     TITLE = {Noncommutative counterparts of the {S}pringer resolution},
 BOOKTITLE = {International {C}ongress of {M}athematicians. {V}ol. {II}},
     PAGES = {1119--1144},
 PUBLISHER = {Eur. Math. Soc., Z\"urich},
      YEAR = {2006},
}

@article{romahecke,
    AUTHOR = {Bezrukavnikov, Roman},
     TITLE = {On two geometric realizations of an affine {H}ecke algebra},
   JOURNAL = {Publ. Math. Inst. Hautes \'{E}tudes Sci.},
  FJOURNAL = {Publications Math\'{e}matiques. Institut de Hautes \'{E}tudes
              Scientifiques},
    VOLUME = {123},
      YEAR = {2016},
     PAGES = {1--67},
     eprint={1209.0403},
}

@article{CaraianiCEB,
     Author={Caraiani, Ana},
     Title = {An excursion through the land of shtukas},
     Note={\url{https://www.ams.org/meetings/lectures/current-events-bulletin}}
}

@book{HA,
	Author = {Lurie, Jacob},
	Title = {Higher algebra},
        Note = {\url{http://math.ias.edu/~lurie/papers/HA.pdf}}
}

@book {HTT,
    AUTHOR = {Lurie, Jacob},
     TITLE = {Higher topos theory},
    SERIES = {Annals of Mathematics Studies},
    VOLUME = {170},
 PUBLISHER = {Princeton University Press, Princeton, NJ},
      YEAR = {2009},
     PAGES = {xviii+925},
     Note = {\url{http://math.ias.edu/~lurie/papers/HTT.pdf}}  
}

@article {MV,
    AUTHOR = {Mirkovi\'{c}, Ivan and Vilonen, Kari},
     TITLE = {Geometric {L}anglands duality and representations of algebraic
              groups over commutative rings},
   JOURNAL = {Ann. of Math. (2)},
  FJOURNAL = {Annals of Mathematics. Second Series},
    VOLUME = {166},
      YEAR = {2007},
    NUMBER = {1},
     PAGES = {95--143},
       DOI = {10.4007/annals.2007.166.95},
     }

@incollection {FeiginFrenkel,
    AUTHOR = {Feigin, Boris and Frenkel, Edward},
     TITLE = {Affine {K}ac-{M}oody algebras at the critical level and
              {G}elfand-{D}iki\u i\ algebras},
 BOOKTITLE = {Infinite analysis, {P}art {A}, {B} ({K}yoto, 1991)},
    SERIES = {Adv. Ser. Math. Phys.},
    VOLUME = {16},
     PAGES = {197--215},
 PUBLISHER = {World Sci. Publ., River Edge, NJ},
      YEAR = {1992},
}

@book {FrenkelLanglands,
    AUTHOR = {Frenkel, Edward},
     TITLE = {Langlands correspondence for loop groups},
    SERIES = {Cambridge Studies in Advanced Mathematics},
    VOLUME = {103},
 PUBLISHER = {Cambridge University Press, Cambridge},
      YEAR = {2007},
     PAGES = {xvi+379},
}

@article{ayalafrancis,
	Author = {David Ayala and John Francis},
	Journal = {J. of Topology},
	Number = {4},
	Pages = {1045-1084},
	Title = {Factorization homology of topological manifolds},
	Volume = {8},
	Year = {2015},
       DOI = {10.1112/jtopol/jtv028},
}

@article{ArinkinGaitsgory,
	Author = {Dima Arinkin and Dennis Gaitsgory},
	Journal = {Selecta Math. (N.S.)},
	Number = {1},
	Pages = {1-199},
	Title = {Singular support of coherent sheaves, and the geometric {Langlands} conjecture.},
	Volume = {21},
	Year = {2015},
       DOI = {10.1007/s00029-014-0167-5},}

@article {Arinkinlambda,
    AUTHOR = {Arinkin, D.},
     TITLE = {On {$\lambda$}-connections on a curve where {$\lambda$} is a
              formal parameter},
   JOURNAL = {Math. Res. Lett.},
  FJOURNAL = {Mathematical Research Letters},
    VOLUME = {12},
      YEAR = {2005},
    NUMBER = {4},
     PAGES = {551--565},}

@article{BezFink,
	Author = {Roman Bezrukavnikov and Mikhail Finkelberg},
	Journal = {Mosc. Math. J. 8 (2008), no. 1, 39--72, 183.},
	Number = {1},
	Pages = {39-72},
	Title = {Equivariant {S}atake category and {K}ostant-{W}hittaker reduction.},
	Volume = {8},
	Year = {2008},
       DOI = {10.17323/1609-4514-2008-8-1-39-72},}

@article{BettiLanglands,
	Author = {David Ben-Zvi and David Nadler},
	Journal = {Proc. Sympos. Pure Math.},
	Title = {Betti Geometric {L}anglands},
	Volume = {97.2},
	Year = {2018},
       DOI = {10.1090/pspum/097.2/01},
}

@article {Rothstein,
    AUTHOR = {Rothstein, Mitchell},
     TITLE = {Sheaves with connection on abelian varieties},
   JOURNAL = {Duke Math. J.},
  FJOURNAL = {Duke Mathematical Journal},
    VOLUME = {84},
      YEAR = {1996},
    NUMBER = {3},
     PAGES = {565--598},
       DOI = {10.1215/S0012-7094-96-08418-5},
}

@misc{LaumonFourier,
      title={Transformation de {F}ourier generalis\'ee}, 
      author={Gerard Laumon},
      year={1996},
      eprint={alg-geom/9603004},
      archivePrefix={arXiv},
      primaryClass={alg-geom}
}

@article {LaumonNilpotent,
    AUTHOR = {Laumon, G\'erard},
     TITLE = {Un analogue global du c\^one nilpotent},
   JOURNAL = {Duke Math. J.},
  FJOURNAL = {Duke Mathematical Journal},
    VOLUME = {57},
      YEAR = {1988},
    NUMBER = {2},
     PAGES = {647--671}
}

@article {LaumonGLC,
    AUTHOR = {Laumon, G\'erard},
     TITLE = {Correspondance de {L}anglands g\'eom\'etrique pour les corps
              de fonctions},
   JOURNAL = {Duke Math. J.},
  FJOURNAL = {Duke Mathematical Journal},
    VOLUME = {54},
      YEAR = {1987},
    NUMBER = {2},
     PAGES = {309--359}
}

@book{CostelloGwilliam,
    AUTHOR = {Costello, Kevin and Gwilliam, Owen},
     TITLE = {Factorization algebras in quantum field theory. {V}ol. 1},
    SERIES = {New Mathematical Monographs},
    VOLUME = {31},
 PUBLISHER = {Cambridge University Press, Cambridge},
      YEAR = {2017},
     PAGES = {ix+387},
      ISBN = {978-1-107-16310-2},
   MRCLASS = {81-01 (17B69 18D50 53D55 81R05 81R10)},
  MRNUMBER = {3586504},
MRREVIEWER = {Domenico\ Fiorenza},
       DOI = {10.1017/9781316678626},
}

@article {Frenkelgauge,
    AUTHOR = {Frenkel, Edward},
     TITLE = {Gauge theory and {L}anglands duality},
      NOTE = {S\'eminaire Bourbaki. Volume 2008/2009. Expos\'es 997--1011},
   JOURNAL = {Ast\'erisque},
  FJOURNAL = {Ast\'erisque},
    NUMBER = {332},
      YEAR = {2010},
     PAGES = {Exp. No. 1010, ix--x, 369--403},
      ISSN = {0303-1179,2492-5926},
      ISBN = {978-2-85629-291-4},
   MRCLASS = {22E57 (14D24 81T60)},
  MRNUMBER = {2648685},
MRREVIEWER = {Carlos\ T.\ Simpson},
}

@article {GaiottoWittenboundary,
    AUTHOR = {Gaiotto, Davide and Witten, Edward},
     TITLE = {Supersymmetric boundary conditions in {$\mathcal N=4$} super
              {Y}ang-{M}ills theory},
   JOURNAL = {J. Stat. Phys.},
  FJOURNAL = {Journal of Statistical Physics},
    VOLUME = {135},
      YEAR = {2009},
    NUMBER = {5-6},
     PAGES = {789--855},
       DOI = {10.1007/s10955-009-9687-3},
}

@article {GaiottoWittenSduality,
    AUTHOR = {Gaiotto, Davide and Witten, Edward},
     TITLE = {{$S$}-duality of boundary conditions in {$\mathcal N=4$} super
              {Y}ang-{M}ills theory},
   JOURNAL = {Adv. Theor. Math. Phys.},
  FJOURNAL = {Advances in Theoretical and Mathematical Physics},
    VOLUME = {13},
      YEAR = {2009},
    NUMBER = {3},
     PAGES = {721--896},
}

@article {FrenkelGaiotto,
    AUTHOR = {Frenkel, Edward and Gaiotto, Davide},
     TITLE = {Quantum {L}anglands dualities of boundary conditions,
              {$D$}-modules, and conformal blocks},
   JOURNAL = {Commun. Number Theory Phys.},
  FJOURNAL = {Communications in Number Theory and Physics},
    VOLUME = {14},
      YEAR = {2020},
    NUMBER = {2},
     PAGES = {199--313},
       DOI = {10.4310/CNTP.2020.v14.n2.a1},
}

@article{Ginzburg1995,	
      Author = {Victor Ginzburg},	
      Eprint = {alg-geom/9511007v4},
	Eprintclass = {alg-geom},
	Eprinttype = {arXiv},
	Title = {Perverse sheaves on a Loop group and Langlands' duality}}

@article{KapustinWitten,
	Archiveprefix = {arXiv},
	Author = {Kapustin, Anton and Witten, Edward},
	Doi = {10.4310/CNTP.2007.v1.n1.a1},
	Journal = {Commun. Num. Theor. Phys.},
	Pages = {1-236},
	Primaryclass = {hep-th},
	Title = {Electric-Magnetic Duality And The Geometric {L}anglands Program},
	Volume = {1},
	Year = {2007}}

@article {HitchinSection,
    AUTHOR = {Hitchin, N. J.},
     TITLE = {Lie groups and {T}eichm\"uller space},
   JOURNAL = {Topology},
  FJOURNAL = {Topology. An International Journal of Mathematics},
    VOLUME = {31},
      YEAR = {1992},
    NUMBER = {3},
     PAGES = {449--473},
}

@incollection {HalesBourbaki,
    AUTHOR = {Hales, Thomas C.},
     TITLE = {The fundamental lemma and the {H}itchin fibration [after
              {N}g\^o{} {B}ao {C}h\^au]},
      NOTE = {S\'eminaire Bourbaki: Vol. 2010/2011. Expos\'es 1027--1042},
   JOURNAL = {Ast\'erisque},
  FJOURNAL = {Ast\'erisque},
    NUMBER = {348},
      YEAR = {2012},
     PAGES = {Exp. No. 1035, ix, 233--263},}

@incollection {StrohBourbaki,
    AUTHOR = {Stroh, Beno\^it},
     TITLE = {La param\'etrisation de {L}anglands globale sur les corps de
              fonctions},
      NOTE = {S\'eminaire Bourbaki. Vol. 2015/2016. Expos\'es 1104--1119},
   JOURNAL = {Ast\'erisque},
  FJOURNAL = {Ast\'erisque},
    NUMBER = {390},
      YEAR = {2017},
     PAGES = {Exp. No. 1110, 169--197},
}

@article{AGKRRV1,
	Archiveprefix = {arXiv},
	Author = {Arinkin, Dima and Gaitsgory, Dennis and Kazhdan, David and Raskin, Sam and Rozenblyum, Nick and Varshavsky, Yakov},
	Eprint = {2010.01906},
	Title = {The stack of local systems with restricted variation and geometric {L}anglands theory with nilpotent singular support},
	Year = {2020}}

@article{AGKRRV3,
	Archiveprefix = {arXiv},
	Author = {Arinkin, Dima and Gaitsgory, Dennis and Kazhdan, David and Raskin, Sam and Rozenblyum, Nick and Varshavsky, Yakov},
	Eprint = {2102.07906},
	Title = {Automorphic functions as the trace of {F}robenius},
	Year = {2021}}

@book {GR,
    AUTHOR = {Gaitsgory, Dennis and Rozenblyum, Nick},
     TITLE = {A study in derived algebraic geometry. {V}ol. {I}.
              {C}orrespondences and duality, {V}ol. {II}.
              {D}eformations, {L}ie theory and formal geometry},
    SERIES = {Mathematical Surveys and Monographs},
    VOLUME = {221},
 PUBLISHER = {American Mathematical Society, Providence, RI},
      YEAR = {2017}
}

@article{FGV,
	Author = {Frenkel, Edward and Gaitsgory, Dennis and Vilonen, Kari},
	Journal={Ann. Math},	
	Title = {Whittaker patterns in the geometry of moduli spaces of bundles on
curves},
	issue={153},
	pages={699--748},
	Year = {2001},
       DOI = {10.2307/2661366},}

@incollection {donagimarkman,
    AUTHOR = {Donagi, Ron and Markman, Eyal},
     TITLE = {Spectral covers, algebraically completely integrable,
              {H}amiltonian systems, and moduli of bundles},
 BOOKTITLE = {Integrable systems and quantum groups ({M}ontecatini {T}erme,
              1993)},
    SERIES = {Lecture Notes in Math.},
    VOLUME = {1620},
     PAGES = {1--119},
 PUBLISHER = {Springer, Berlin},
      YEAR = {1996}
}

@article {donagigaitsgory,
    AUTHOR = {Donagi, Ron Y. and Gaitsgory, Dennis},
     TITLE = {The gerbe of {H}iggs bundles},
   JOURNAL = {Transform. Groups},
  FJOURNAL = {Transformation Groups},
    VOLUME = {7},
      YEAR = {2002},
    NUMBER = {2},
     PAGES = {109--153},
}

@article {Mukai,
    AUTHOR = {Mukai, Shigeru},
     TITLE = {Duality between {$D(X)$}\ and {$D(\hat X)$}\ with its
              application to {P}icard sheaves},
   JOURNAL = {Nagoya Math. J.},
  FJOURNAL = {Nagoya Mathematical Journal},
    VOLUME = {81},
      YEAR = {1981},
     PAGES = {153--175},
      ISSN = {0027-7630,2152-6842},
   MRCLASS = {14K05 (14H40)},
  MRNUMBER = {607081},
MRREVIEWER = {Allen\ B.\ Altman},
       URL = {http://projecteuclid.org.ezproxy.lib.utexas.edu/euclid.nmj/1118786312},
}

@article {donagipantev,
    AUTHOR = {Donagi, Ron and Pantev, Tony},
     TITLE = {Langlands duality for {H}itchin systems},
   JOURNAL = {Invent. Math.},
  FJOURNAL = {Inventiones Mathematicae},
    VOLUME = {189},
      YEAR = {2012},
    NUMBER = {3},
     PAGES = {653--735},}

@article {arinkinFM,
    AUTHOR = {Arinkin, Dima},
     TITLE = {Autoduality of compactified {J}acobians for curves with plane
              singularities},
   JOURNAL = {J. Algebraic Geom.},
  FJOURNAL = {Journal of Algebraic Geometry},
    VOLUME = {22},
      YEAR = {2013},
    NUMBER = {2},
     PAGES = {363--388},
      ISSN = {1056-3911,1534-7486},
   MRCLASS = {14H40 (14D20 14F05 14H20 14K30)},
  MRNUMBER = {3019453},
MRREVIEWER = {H.\ Lange},
       DOI = {10.1090/S1056-3911-2012-00596-7},
       URL = {https://doi-org.ezproxy.lib.utexas.edu/10.1090/S1056-3911-2012-00596-7},
}

@article {MRV2,
    AUTHOR = {Melo, Margarida and Rapagnetta, Antonio and Viviani, Filippo},
     TITLE = {Fourier-{M}ukai and autoduality for compactified {J}acobians,
              {II}},
   JOURNAL = {Geom. Topol.},
  FJOURNAL = {Geometry \& Topology},
    VOLUME = {23},
      YEAR = {2019},
    NUMBER = {5},
     PAGES = {2335--2395}  
}

@article {MRV1,
    AUTHOR = {Melo, Margarida and Rapagnetta, Antonio and Viviani, Filippo},
     TITLE = {Fourier-{M}ukai and autoduality for compactified {J}acobians.
              {I}},
   JOURNAL = {J. Reine Angew. Math.},
  FJOURNAL = {Journal f\"ur die Reine und Angewandte Mathematik. [Crelle's
              Journal]},
    VOLUME = {755},
      YEAR = {2019},
     PAGES = {1--65}}

@misc{BZSV,
      title={Relative Langlands Duality}, 
      author={David Ben-Zvi and Yiannis Sakellaridis and Akshay Venkatesh},
      year={2024},
      eprint={2409.04677},
      archivePrefix={arXiv},
      primaryClass={math.RT},
      url={https://arxiv.org/abs/2409.04677}, 
}

@article {HauselThaddeus,
    AUTHOR = {Hausel, Tam\'as and Thaddeus, Michael},
     TITLE = {Mirror symmetry, {L}anglands duality, and the {H}itchin
              system},
   JOURNAL = {Invent. Math.},
  FJOURNAL = {Inventiones Mathematicae},
    VOLUME = {153},
      YEAR = {2003},
    NUMBER = {1},
     PAGES = {197--229},
      ISSN = {0020-9910,1432-1297},
   MRCLASS = {14J32 (14D20 32J25 32Q25)},
  MRNUMBER = {1990670},
MRREVIEWER = {Anatoly\ Libgober},
       DOI = {10.1007/s00222-003-0286-7},
       URL = {https://doi-org.ezproxy.lib.utexas.edu/10.1007/s00222-003-0286-7},
}

@incollection {HauselICM,
    AUTHOR = {Hausel, Tam\'as},
     TITLE = {Enhanced mirror symmetry for {L}anglands dual {H}itchin
              systems},
 BOOKTITLE = {I{CM}---{I}nternational {C}ongress of {M}athematicians. {V}ol.
              3. {S}ections 1--4},
     PAGES = {2228--2249},
 PUBLISHER = {EMS Press, Berlin},
      YEAR = {[2023] \copyright 2023},
      ISBN = {978-3-98547-061-7; 978-3-98547-561-2; 978-3-98547-058-7},
   MRCLASS = {14D21 (14C17 14D24 14J33 20G05 81T30)},
  MRNUMBER = {4680316},
}

@incollection {EFK,
    AUTHOR = {Etingof, Pavel and Frenkel, Edward and Kazhdan, David},
     TITLE = {An analytic version of the {L}anglands correspondence for
              complex curves},
 BOOKTITLE = {Integrability, quantization, and geometry {II}. {Q}uantum
              theories and algebraic geometry},
    SERIES = {Proc. Sympos. Pure Math.},
    VOLUME = {103.2},
     PAGES = {137--202},
 PUBLISHER = {Amer. Math. Soc., Providence, RI},
      YEAR = {[2021] \copyright 2021},
      ISBN = {978-1-4704-5592-7},
   MRCLASS = {14D24 (11R39 22E57 35P05)},
  MRNUMBER = {4285697},
       DOI = {10.1090/pspum/103.2/01856},
       URL = {https://doi-org.ezproxy.lib.utexas.edu/10.1090/pspum/103.2/01856},
}

@misc{jordansurvey,
      title={Langlands duality for skein modules of 3-manifolds}, 
      author={David Jordan},
      year={2023},
      eprint={2302.14734},
      archivePrefix={arXiv},
      primaryClass={math.QA},
      url={https://arxiv.org/abs/2302.14734}, 
}

@article {GaiottoWittenAnalytic,
    AUTHOR = {Gaiotto, Davide and Witten, Edward},
     TITLE = {Gauge theory and the analytic form of the geometric
              {L}anglands program},
   JOURNAL = {Ann. Henri Poincar\'e},
  FJOURNAL = {Annales Henri Poincar\'e. A Journal of Theoretical and
              Mathematical Physics},
    VOLUME = {25},
      YEAR = {2024},
    NUMBER = {1},
     PAGES = {557--671},
      ISSN = {1424-0637,1424-0661},
   MRCLASS = {14D24 (53D30 81T13 81T30)},
  MRNUMBER = {4699908},
MRREVIEWER = {K.\ Maruyoshi},
       DOI = {10.1007/s00023-022-01225-6},
       URL = {https://doi-org.ezproxy.lib.utexas.edu/10.1007/s00023-022-01225-6},
}

@article {AtiyahBottYM,
    AUTHOR = {Atiyah, M. F. and Bott, R.},
     TITLE = {The {Y}ang-{M}ills equations over {R}iemann surfaces},
   JOURNAL = {Philos. Trans. Roy. Soc. London Ser. A},
  FJOURNAL = {Philosophical Transactions of the Royal Society of London.
              Series A. Mathematical and Physical Sciences},
    VOLUME = {308},
      YEAR = {1983},
    NUMBER = {1505},
     PAGES = {523--615}
}

@misc{scholzemotivic,
      title={Geometrization of the local Langlands correspondence, motivically}, 
      author={Peter Scholze},
      year={2025},
      eprint={2501.07944},
      archivePrefix={arXiv},
      primaryClass={math.AG},
      url={https://arxiv.org/abs/2501.07944}, 
}

@incollection {donagipantevsurvey,
    AUTHOR = {Donagi, R. and Pantev, T.},
     TITLE = {Geometric {L}anglands and non-abelian {H}odge theory},
 BOOKTITLE = {Surveys in differential geometry. {V}ol. {XIII}. {G}eometry,
              analysis, and algebraic geometry: forty years of the {J}ournal
              of {D}ifferential {G}eometry},
    SERIES = {Surv. Differ. Geom.},
    VOLUME = {13},
     PAGES = {85--116},
 PUBLISHER = {Int. Press, Somerville, MA},
      YEAR = {2009},}

@misc{richesurvey,
      title={Some applications of the geometric Satake equivalence to modular representation theory}, 
      author={Simon Riche},
      year={2024},
      eprint={2403.03734},
      archivePrefix={arXiv},
      primaryClass={math.RT},
      url={https://arxiv.org/abs/2403.03734}, 
}

@misc{GLCI,
      title={Proof of the geometric Langlands conjecture I: construction of the functor}, 
      author={Dennis Gaitsgory and Sam Raskin},
      year={2024},
      eprint={2405.03599},
      archivePrefix={arXiv},
      primaryClass={math.AG},
      url={https://arxiv.org/abs/2405.03599}, 
}

@misc{GLCII,
      title={Proof of the geometric Langlands conjecture II: Kac-Moody localization and the FLE}, 
      author={Dima Arinkin and Dario Beraldo and Justin Campbell and Lin Chen and Joakim Faergeman and Dennis Gaitsgory and Kevin Lin and Sam Raskin and Nick Rozenblyum},
      year={2024},
      eprint={2409.07051},
      archivePrefix={arXiv},
      primaryClass={math.AG},
      url={https://arxiv.org/abs/2409.07051}, 
}

@misc{GLCIII,
      title={Proof of the geometric Langlands conjecture III: compatibility with parabolic induction}, 
      author={Justin Campbell and Lin Chen and Joakim Faergeman and Dennis Gaitsgory and Kevin Lin and Sam Raskin and Nick Rozenblyum},
      year={2024},
      eprint={2409.07051},
      archivePrefix={arXiv},
      primaryClass={math.AG},
      url={https://arxiv.org/abs/2409.07051}, 
}

@misc{GLCIV,
      title={Proof of the geometric Langlands conjecture IV: ambidexterity}, 
      author={D. Arinkin and D. Beraldo and L. Chen and J. Faergeman and D. Gaitsgory and K. Lin and S. Raskin and N. Rozenblyum},
      year={2024},
      eprint={2409.08670},
      archivePrefix={arXiv},
      primaryClass={math.AG},
      url={https://arxiv.org/abs/2409.08670}, 
}

@misc{GLCV,
      title={Proof of the geometric Langlands conjecture V: the multiplicity one theorem}, 
      author={Dennis Gaitsgory and Sam Raskin},
      year={2024},
      eprint={2409.09856},
      archivePrefix={arXiv},
      primaryClass={math.AG},
      url={https://arxiv.org/abs/2409.09856}, 
}

@misc{gaitsgoryraskin2025,
      title={Geometric Langlands in positive characteristic from characteristic zero}, 
      author={Dennis Gaitsgory and Sam Raskin},
      year={2025},
      eprint={2508.02237},
      archivePrefix={arXiv},
      primaryClass={math.AG},
      url={https://arxiv.org/abs/2508.02237}, 
}

@book{scholzetwistor,
	Author = {Scholze, Peter},
	Title = {Geometrization of the real local Langlands
correspondence},
        Note = {\url{https://people.mpim-bonn.mpg.de/scholze/RealLocalLanglands.pdf}},
}

@article {NadlerCEB,
    AUTHOR = {Nadler, David},
     TITLE = {The geometric nature of the fundamental lemma},
   JOURNAL = {Bull. Amer. Math. Soc. (N.S.)},
  FJOURNAL = {American Mathematical Society. Bulletin. New Series},
    VOLUME = {49},
      YEAR = {2012},
    NUMBER = {1},
     PAGES = {1--50},
      ISSN = {0273-0979,1088-9485},
   MRCLASS = {11R39 (11F72 14D24 22E35)},
  MRNUMBER = {2869006},
MRREVIEWER = {Anne-Marie\ H.\ Aubert},
       DOI = {10.1090/S0273-0979-2011-01342-8},
       URL = {https://doi.org/10.1090/S0273-0979-2011-01342-8},
}

@article {FrenkelCEB,
    AUTHOR = {Frenkel, Edward},
     TITLE = {Recent advances in the {L}anglands program},
   JOURNAL = {Bull. Amer. Math. Soc. (N.S.)},
  FJOURNAL = {American Mathematical Society. Bulletin. New Series},
    VOLUME = {41},
      YEAR = {2004},
    NUMBER = {2},
     PAGES = {151--184},
      ISSN = {0273-0979,1088-9485},
   MRCLASS = {11R39 (11F52 11G45 14D20)},
  MRNUMBER = {2043750},
MRREVIEWER = {Igor\ Yu.\ Potemine},
       DOI = {10.1090/S0273-0979-04-01001-8},
       URL = {https://doi.org/10.1090/S0273-0979-04-01001-8},
}

@article {NadlerZaslow,
    AUTHOR = {Nadler, David and Zaslow, Eric},
     TITLE = {Constructible sheaves and the {F}ukaya category},
   JOURNAL = {J. Amer. Math. Soc.},
  FJOURNAL = {Journal of the American Mathematical Society},
    VOLUME = {22},
      YEAR = {2009},
    NUMBER = {1},
     PAGES = {233--286},
      ISSN = {0894-0347,1088-6834},
   MRCLASS = {53D37 (32S60 53D40)},
  MRNUMBER = {2449059},
MRREVIEWER = {Richard\ P.\ Thomas},
       DOI = {10.1090/S0894-0347-08-00612-7},
       URL = {https://doi.org/10.1090/S0894-0347-08-00612-7},
}

@article {NadlerBrane,
    AUTHOR = {Nadler, David},
     TITLE = {Microlocal branes are constructible sheaves},
   JOURNAL = {Selecta Math. (N.S.)},
  FJOURNAL = {Selecta Mathematica. New Series},
    VOLUME = {15},
      YEAR = {2009},
    NUMBER = {4},
     PAGES = {563--619}
}

@article {GPS,
    AUTHOR = {Ganatra, Sheel and Pardon, John and Shende, Vivek},
     TITLE = {Microlocal {M}orse theory of wrapped {F}ukaya categories},
   JOURNAL = {Ann. of Math. (2)},
  FJOURNAL = {Annals of Mathematics. Second Series},
    VOLUME = {199},
      YEAR = {2024},
    NUMBER = {3},
     PAGES = {943--1042}
}

@article{Kontsevich,
	Author = {Kontsevich, Maxim},
	Title = {Symplectic geometry of homological algebra},
	Year = {2009},
        URL = {https://archive.mpim-bonn.mpg.de/id/eprint/1536/1/preprint_2009_40_a.pdf},
}

@article {FR,
    AUTHOR = {F\ae rgeman, Joakim and Raskin, Sam},
     TITLE = {Non-vanishing of geometric {W}hittaker coefficients for
              reductive groups},
   JOURNAL = {J. Amer. Math. Soc.},
  FJOURNAL = {Journal of the American Mathematical Society},
    VOLUME = {38},
      YEAR = {2025},
    NUMBER = {4},
     PAGES = {919--995},
}

@incollection {GaitsgoryBourbaki,
    AUTHOR = {Gaitsgory, Dennis},
     TITLE = {Progr\`es r\'ecents dans la th\'eorie de {L}anglands
              g\'eom\'etrique},
      NOTE = {S\'eminaire Bourbaki. Vol. 2015/2016. Expos\'es 1104--1119},
   JOURNAL = {Ast\'erisque},
  FJOURNAL = {Ast\'erisque},
    NUMBER = {390},
      YEAR = {2017},
     PAGES = {Exp. No. 1109, 139--168},
      ISSN = {0303-1179,2492-5926},
      ISBN = {978-2-85629-855-8},
   MRCLASS = {14D24 (14H60)},
  MRNUMBER = {3666025},
MRREVIEWER = {Andrei\ D.\ Halanay},
}

@misc{BKS,
      title={Contractibility of the space of generic opers for classical groups}, 
      author={Dario Beraldo and David Kazhdan and Tomer M. Schlank},
      year={2022},
      eprint={1801.00655},
      archivePrefix={arXiv},
      primaryClass={math.RT},
      url={https://arxiv.org/abs/1801.00655}, 
}

@article {zhiweiGalois,
    AUTHOR = {Yun, Zhiwei},
     TITLE = {Motives with exceptional {G}alois groups and the inverse
              {G}alois problem},
   JOURNAL = {Invent. Math.},
  FJOURNAL = {Inventiones Mathematicae},
    VOLUME = {196},
      YEAR = {2014},
    NUMBER = {2},
     PAGES = {267--337},
      ISSN = {0020-9910,1432-1297},
   MRCLASS = {14D24 (12F12 14C15 14G32 20G41)},
  MRNUMBER = {3193750},
MRREVIEWER = {Satoshi\ Mochizuki},
       DOI = {10.1007/s00222-013-0469-9},
       URL = {https://doi.org/10.1007/s00222-013-0469-9},
}

@inproceedings {zhiweirigidity,
    AUTHOR = {Yun, Zhiwei},
     TITLE = {Rigidity in the {L}anglands correspondence and applications},
 BOOKTITLE = {Proceedings of the {S}ixth {I}nternational {C}ongress of
              {C}hinese {M}athematicians. {V}ol. {I}},
    SERIES = {Adv. Lect. Math. (ALM)},
    VOLUME = {36},
     PAGES = {199--234},
 PUBLISHER = {Int. Press, Somerville, MA},
      YEAR = {2017},
      ISBN = {978-1-57146-348-7},
   MRCLASS = {14D24 (11S37 12F12 14H25 20G41)},
  MRNUMBER = {3702076},
MRREVIEWER = {Volker\ J.\ Heiermann},
}

@article {LamTemplier,
    AUTHOR = {Lam, Thomas and Templier, Nicolas},
     TITLE = {The mirror conjecture for minuscule flag varieties},
   JOURNAL = {Duke Math. J.},
  FJOURNAL = {Duke Mathematical Journal},
    VOLUME = {173},
      YEAR = {2024},
    NUMBER = {1},
     PAGES = {75--175},
      ISSN = {0012-7094,1547-7398},
   MRCLASS = {14D24 (11T23 14J33 14M15)},
  MRNUMBER = {4728689},
       DOI = {10.1215/00127094-2024-0007},
       URL = {https://doi.org/10.1215/00127094-2024-0007},
}

@article {FrenkelGross,
    AUTHOR = {Frenkel, Edward and Gross, Benedict},
     TITLE = {A rigid irregular connection on the projective line},
   JOURNAL = {Ann. of Math. (2)},
  FJOURNAL = {Annals of Mathematics. Second Series},
    VOLUME = {170},
      YEAR = {2009},
    NUMBER = {3},
     PAGES = {1469--1512},
      ISSN = {0003-486X,1939-8980},
   MRCLASS = {14D24 (14F40)},
  MRNUMBER = {2600880},
       DOI = {10.4007/annals.2009.170.1469},
       URL = {https://doi.org/10.4007/annals.2009.170.1469},
}

@article {HeinlochNgoYun,
    AUTHOR = {Heinloth, Jochen and Ng\^o, Bao-Ch\^au and Yun, Zhiwei},
     TITLE = {Kloosterman sheaves for reductive groups},
   JOURNAL = {Ann. of Math. (2)},
  FJOURNAL = {Annals of Mathematics. Second Series},
    VOLUME = {177},
      YEAR = {2013},
    NUMBER = {1},
     PAGES = {241--310},
      ISSN = {0003-486X,1939-8980},
   MRCLASS = {22E57 (11L05)},
  MRNUMBER = {2999041},
MRREVIEWER = {K.\ Strambach},
       DOI = {10.4007/annals.2013.177.1.5},
       URL = {https://doi.org/10.4007/annals.2013.177.1.5},
}

@incollection {zhuICM,
    AUTHOR = {Zhu, Xinwen},
     TITLE = {Arithmetic and geometric {L}anglands program},
 BOOKTITLE = {I{CM}---{I}nternational {C}ongress of {M}athematicians. {V}ol.
              3. {S}ections 1--4},
     PAGES = {2012--2045},
 PUBLISHER = {EMS Press, Berlin},
      YEAR = {[2023] \copyright 2023},
      ISBN = {978-3-98547-061-7; 978-3-98547-561-2; 978-3-98547-058-7},
   MRCLASS = {11R39 (11G18 11G40 11S37 14D24)},
  MRNUMBER = {4680309},
MRREVIEWER = {Feng-Wen\ An},
}

@incollection {FrenkelGaudin,
    AUTHOR = {Frenkel, Edward},
     TITLE = {Gaudin model and opers},
 BOOKTITLE = {Infinite dimensional algebras and quantum integrable systems},
    SERIES = {Progr. Math.},
    VOLUME = {237},
     PAGES = {1--58},
 PUBLISHER = {Birkh\"auser, Basel},
      YEAR = {2005},
      ISBN = {978-3-7643-7215-6; 3-7643-7215-X},
   MRCLASS = {17B67 (37K30 81R12 82-02 82B23)},
  MRNUMBER = {2160841},
       DOI = {10.1007/3-7643-7341-5\_1},
       URL = {https://doi.org/10.1007/3-7643-7341-5_1},
}

@article {XinwenBessel,
    AUTHOR = {Xu, Daxin and Zhu, Xinwen},
     TITLE = {Bessel {$F$}-isocrystals for reductive groups},
   JOURNAL = {Invent. Math.},
  FJOURNAL = {Inventiones Mathematicae},
    VOLUME = {227},
      YEAR = {2022},
    NUMBER = {3},
     PAGES = {997--1092},
      ISSN = {0020-9910,1432-1297},
   MRCLASS = {14F30 (14D24 14F10)},
  MRNUMBER = {4384193},
MRREVIEWER = {Nguy\cftil en Qu\^oc Th\'ang},
       DOI = {10.1007/s00222-021-01079-5},
       URL = {https://doi.org/10.1007/s00222-021-01079-5},
}

@article {LiNadlerYun,
    AUTHOR = {Li, Penghui and Nadler, David and Yun, Zhiwei},
     TITLE = {Functions on the commuting stack via {L}anglands duality},
   JOURNAL = {Ann. of Math. (2)},
  FJOURNAL = {Annals of Mathematics. Second Series},
    VOLUME = {200},
      YEAR = {2024},
    NUMBER = {2},
     PAGES = {609--748},
      ISSN = {0003-486X,1939-8980},
   MRCLASS = {14D24 (20G99)},
  MRNUMBER = {4792071},
       DOI = {10.4007/annals.2024.200.2.5},
       URL = {https://doi.org/10.4007/annals.2024.200.2.5},
}

@article {CarterLusztig,
    AUTHOR = {Carter, R. W.},
     TITLE = {A survey of the work of {G}eorge {L}usztig},
   JOURNAL = {Nagoya Math. J.},
  FJOURNAL = {Nagoya Mathematical Journal},
    VOLUME = {182},
      YEAR = {2006},
     PAGES = {1--45},
}

@article {BezBrav,
    AUTHOR = {Braverman, Alexander and Bezrukavnikov, Roman},
     TITLE = {Geometric {L}anglands correspondence for {$\mathcal D$}-modules in
              prime characteristic: the {${\rm GL}(n)$} case},
   JOURNAL = {Pure Appl. Math. Q.},
  FJOURNAL = {Pure and Applied Mathematics Quarterly},
    VOLUME = {3},
      YEAR = {2007},
    NUMBER = {1},
     PAGES = {153--179},
}

@article {ChenZhu,
    AUTHOR = {Chen, Tsao-Hsien and Zhu, Xinwen},
     TITLE = {Geometric {L}anglands in prime characteristic},
   JOURNAL = {Compos. Math.},
  FJOURNAL = {Compositio Mathematica},
    VOLUME = {153},
      YEAR = {2017},
    NUMBER = {2},
     PAGES = {395--452},
}

@article {DonagiPanteveigensheaves,
    AUTHOR = {Donagi, Ron and Pantev, Tony},
     TITLE = {Parabolic {H}ecke eigensheaves},
   JOURNAL = {Ast\'erisque},
  FJOURNAL = {Ast\'erisque},
    NUMBER = {435},
      YEAR = {2022},
     PAGES = {viii+192},
     }

@misc{DonagiPantevSimpson,
      title={Twistor Hecke eigensheaves in genus 2}, 
      author={Ron Donagi and Tony Pantev and Carlos Simpson},
      year={2024},
      eprint={2403.17045},
      archivePrefix={arXiv},
      primaryClass={math.AG},
      url={https://arxiv.org/abs/2403.17045}, 
}

@misc{arinkinoper,
      title={Irreducible connections admit generic oper structures}, 
      author={Dima Arinkin},
      year={2016},
      eprint={1602.08989},
      archivePrefix={arXiv},
      primaryClass={math.AG},
      url={https://arxiv.org/abs/1602.08989}, 
}

@misc{gaitsgoryICM,
      title={Local and global Langlands conjecture(s) over function fields. To appear, Proceedings of the 2026 ICM}, 
      author={Dennis Gaitsgory},
      year={2025},
      eprint={2509.24902},
      archivePrefix={arXiv},
      primaryClass={math.AG},
      url={https://arxiv.org/abs/2509.24902}, 
}

@misc{raskinICM,
      title={Conjectures of Arthur and Ramanujan for unramified automorphic forms: Announcement. To appear, Proceedings of the 2026 ICM}, 
      author={Sam Raskin},
     }

@article {DrinfeldGLC,
    AUTHOR = {Drinfeld, V. G.},
     TITLE = {Two-dimensional {$l$}-adic representations of the fundamental
              group of a curve over a finite field and automorphic forms on
              {${\rm GL}(2)$}},
   JOURNAL = {Amer. J. Math.},
  FJOURNAL = {American Journal of Mathematics},
    VOLUME = {105},
      YEAR = {1983},
    NUMBER = {1},
     PAGES = {85--114}
}

@article {Gaitsgoryvanishing,
    AUTHOR = {Gaitsgory, D.},
     TITLE = {On a vanishing conjecture appearing in the geometric
              {L}anglands correspondence},
   JOURNAL = {Ann. of Math. (2)},
  FJOURNAL = {Annals of Mathematics. Second Series},
    VOLUME = {160},
      YEAR = {2004},
    NUMBER = {2},
     PAGES = {617--682},
      ISSN = {0003-486X,1939-8980},
   MRCLASS = {11R39 (11F70 14D20 22E55)},
  MRNUMBER = {2123934},
MRREVIEWER = {Peter\ Fiebig},
       DOI = {10.4007/annals.2004.160.617},
       URL = {https://doi-org.ezproxy.lib.utexas.edu/10.4007/annals.2004.160.617},
}

@article {NadlerTaylor,
    AUTHOR = {Nadler, David and Taylor, Jeremy},
     TITLE = {The {W}hittaker {F}unctional {I}s a {S}hifted {M}icrostalk},
   JOURNAL = {Transform. Groups},
  FJOURNAL = {Transformation Groups},
    VOLUME = {30},
      YEAR = {2025},
    NUMBER = {3},
     PAGES = {1425--1450}
}

@misc{PadurariuToda,
      title={The Dolbeault geometric Langlands conjecture via limit categories}, 
      author={Tudor Pădurariu and Yukinobu Toda},
      year={2025},
      eprint={2508.19624},
      archivePrefix={arXiv},
      primaryClass={math.AG},
      url={https://arxiv.org/abs/2508.19624}, 
}

@article {BMR,
    AUTHOR = {Bezrukavnikov, Roman and Mirkovi\'c, Ivan and Rumynin,
              Dmitriy},
     TITLE = {Localization of modules for a semisimple {L}ie algebra in
              prime characteristic},
      NOTE = {With an appendix by Bezrukavnikov and Simon Riche},
   JOURNAL = {Ann. of Math. (2)},
  FJOURNAL = {Annals of Mathematics. Second Series},
    VOLUME = {167},
      YEAR = {2008},
    NUMBER = {3},
     PAGES = {945--991}}

@article {BezMirk,
    AUTHOR = {Bezrukavnikov, Roman and Mirkovi\'c, Ivan},
     TITLE = {Representations of semisimple {L}ie algebras in prime
              characteristic and the noncommutative {S}pringer resolution},
      NOTE = {With an appendix by Eric Sommers},
   JOURNAL = {Ann. of Math. (2)},
  FJOURNAL = {Annals of Mathematics. Second Series},
    VOLUME = {178},
      YEAR = {2013},
    NUMBER = {3},
     PAGES = {835--919}
}

@article {RicheGeordielinkage,
    AUTHOR = {Riche, Simon and Williamson, Geordie},
     TITLE = {Smith-{T}reumann theory and the linkage principle},
   JOURNAL = {Publ. Math. Inst. Hautes \'Etudes Sci.},
  FJOURNAL = {Publications Math\'ematiques. Institut de Hautes \'Etudes
              Scientifiques},
    VOLUME = {136},
      YEAR = {2022},
     PAGES = {225--292}
}

@article {AcharRiche,
    AUTHOR = {Achar, Pramod N. and Riche, Simon},
     TITLE = {Reductive groups, the loop {G}rassmannian, and the {S}pringer
              resolution},
   JOURNAL = {Invent. Math.},
  FJOURNAL = {Inventiones Mathematicae},
    VOLUME = {214},
      YEAR = {2018},
    NUMBER = {1},
     PAGES = {289--436}
}

@article {Geordiesurvey,
    AUTHOR = {Williamson, Geordie},
     TITLE = {Algebraic representations and constructible sheaves},
   JOURNAL = {Jpn. J. Math.},
  FJOURNAL = {Japanese Journal of Mathematics},
    VOLUME = {12},
      YEAR = {2017},
    NUMBER = {2},
     PAGES = {211--259}
}

@article {RaskinYang,
    AUTHOR = {Raskin, Sam and Yang, David},
     TITLE = {Affine {B}eilinson-{B}ernstein localization at the critical
              level},
   JOURNAL = {Ann. of Math. (2)},
  FJOURNAL = {Annals of Mathematics. Second Series},
    VOLUME = {200},
      YEAR = {2024},
    NUMBER = {2},
     PAGES = {487--527}
}

@article {FrenkelGaitsgorylocalization,
    AUTHOR = {Frenkel, Edward and Gaitsgory, Dennis},
     TITLE = {Localization of {$\mathfrak g$}-modules on the affine
              {G}rassmannian},
   JOURNAL = {Ann. of Math. (2)},
  FJOURNAL = {Annals of Mathematics. Second Series},
    VOLUME = {170},
      YEAR = {2009},
    NUMBER = {3},
     PAGES = {1339--1381}}

@incollection {FrenkelGaitsgoryspherical,
    AUTHOR = {Frenkel, Edward and Gaitsgory, Dennis},
     TITLE = {Local geometric {L}anglands correspondence: the spherical
              case},
 BOOKTITLE = {Algebraic analysis and around},
    SERIES = {Adv. Stud. Pure Math.},
    VOLUME = {54},
     PAGES = {167--186},
 PUBLISHER = {Math. Soc. Japan, Tokyo},
      YEAR = {2009},
      ISBN = {978-4-931469-51-8},
   MRCLASS = {22E57},
  MRNUMBER = {2499556},
MRREVIEWER = {Peter\ Fiebig},
       DOI = {10.2969/aspm/05410167},
       URL = {https://doi-org.ezproxy.lib.utexas.edu/10.2969/aspm/05410167},
}

@book {IHESLanglands,
     TITLE = {The {L}anglands {P}rogram},
    SERIES = {Proceedings of Symposia in Pure Mathematics},
    VOLUME = {112.2},
    EDITOR = {Chaudouard, Pierre-Henri and Gan, Wee and Kaletha, Tasho and
              Sakellaridis, Yiannis},
 PUBLISHER = {American Mathematical Society, Providence, RI},
      YEAR = {2025}
}

@book{ganBZSV,
	Author = {Gan, Wee Teck},
	Title = {Relative Langlands Duality, after Ben-Zvi, Sakellaridis and Venkatesh},
	Note={Seminaire Bourbaki 29 Mars 2025, \url{https://www.bourbaki.fr/seminaires/2024-25/Prog_mar-25.html}},
}

@article {FLEtrio,
    AUTHOR = {Campbell, Justin and Dhillon, Gurbir and Raskin, Sam},
     TITLE = {Fundamental local equivalences in quantum geometric
              {L}anglands},
   JOURNAL = {Compos. Math.},
  FJOURNAL = {Compositio Mathematica},
    VOLUME = {157},
      YEAR = {2021},
    NUMBER = {12},
     PAGES = {2699--2732}
}

@misc{gaitsgorylysenkoFLE,
      title={Metaplectic Whittaker category and quantum groups : the "small" FLE}, 
      author={D. Gaitsgory and S. Lysenko},
      year={2019},
      eprint={1903.02279},
      archivePrefix={arXiv},
      primaryClass={math.AG},
      url={https://arxiv.org/abs/1903.02279}, 
}

@article {gaitsgoryWhittaker,
    AUTHOR = {Gaitsgory, D.},
     TITLE = {Twisted {W}hittaker model and factorizable sheaves},
   JOURNAL = {Selecta Math. (N.S.)},
  FJOURNAL = {Selecta Mathematica. New Series},
    VOLUME = {13},
      YEAR = {2008},
    NUMBER = {4},
     PAGES = {617--659},
}

@misc{chenfu,
      title={An Extension of the Kazhdan-Lusztig Equivalence}, 
      author={Lin Chen and Yuchen Fu},
      year={2021},
      eprint={2111.14606},
      archivePrefix={arXiv},
      primaryClass={math.RT},
      url={https://arxiv.org/abs/2111.14606}, 
}

@article {KazhdanLusztig,
    AUTHOR = {Kazhdan, D. and Lusztig, G.},
     TITLE = {Tensor structures arising from affine {L}ie algebras. {I},
              {II}},
   JOURNAL = {J. Amer. Math. Soc.},
  FJOURNAL = {Journal of the American Mathematical Society},
    VOLUME = {6},
      YEAR = {1993},
    NUMBER = {4}}

@article {FockGoncharov,
    AUTHOR = {Fock, Vladimir and Goncharov, Alexander},
     TITLE = {Moduli spaces of local systems and higher {T}eichm\"uller
              theory},
   JOURNAL = {Publ. Math. Inst. Hautes \'Etudes Sci.},
  FJOURNAL = {Publications Math\'ematiques. Institut de Hautes \'Etudes
              Scientifiques},
    NUMBER = {103},
      YEAR = {2006},
     PAGES = {1--211}
}

@article {GMN,
    AUTHOR = {Gaiotto, Davide and Moore, Gregory W. and Neitzke, Andrew},
     TITLE = {Wall-crossing, {H}itchin systems, and the {WKB} approximation},
   JOURNAL = {Adv. Math.},
  FJOURNAL = {Advances in Mathematics},
    VOLUME = {234},
      YEAR = {2013},
     PAGES = {239--403}
}

\end{document}